%% file: veg-r1.tex
\documentclass[a4paper,11pt]{article}
\usepackage[top=25mm,left=22mm]{geometry}\textwidth15.5cm\textheight23.2cm
\usepackage{amsmath,amssymb,amsthm, graphicx}
\usepackage[utf8]{inputenc}
\usepackage{url}
\usepackage{pstricks}
\usepackage{float,paralist}
\usepackage{tabularx}

\input{./mindef.tex}
\def\p2p{{\tt pde2path}}

\def\ddt{\frac{{\rm d}}{{\rm d}t}}
\def\PMAXP{Pontryagin's Maximum Principle}
\newtheorem{theorem}{Theorem}[section]
\newtheorem{remark}[theorem]{Remark}

\def\brem{\begin{remark}}\def\erem{\end{remark}}
\def\eex{\hfill\mbox{$\rfloor$}}
\renewcommand{\arraystretch}{1.1}\renewcommand{\baselinestretch}{1.15}

\begin{document}
\newcolumntype{C}[1]{>{\centering\arraybackslash}p{#1}}
\title{Optimal harvesting and spatial patterns in a semi arid 
vegetation system 
}
\author{Hannes Uecker\\ \small
Institut f\"ur Mathematik, Universit\"at Oldenburg,\\ 
\small D26111 Oldenburg, hannes.uecker@uni-oldenburg.de}

\normalsize
\maketitle
\begin{abstract} We consider an infinite time horizon spatially distributed 
optimal harvesting 
problem for a vegetation and soil water reaction diffusion system, 
with rainfall as the main external parameter. 
By Pontryagin's maximum principle we derive 
the associated four component canonical system, and numerically 
analyze this and hence 
the optimal control problem in two steps. 
First we numerically compute a rather rich bifurcation structure of 
{\em flat} (spatially homogeneous) and {\em patterned} canonical 
steady states 
(FCSS and PCSS, respectively), in 1D and 2D. Then we compute 
time dependent solutions of the canonical system that connect to 
some FCSS or PCSS. The method is efficient in dealing with 
non-unique canonical steady states, and thus also with 
multiple local maxima of the objective function. 
It turns out that over wide parameter regimes the FCSS, i.e., spatially uniform 
harvesting, are not optimal. Instead, 
controlling the system to a PCSS yields a higher profit. 
Moreover, compared to (a simple model of) private optimization, the social 
control gives a higher yield, and vegetation survives 
for much lower rainfall. {In addition, the computation of the optimal (social) 
control gives an optimal tax to incorporate into the private optimization. }
\end{abstract}
\noi

\section{Introduction and main results}
Vegetation patterns such as spots and stripes appear in ecosystems all over the 
world, in particular in so called semi arid areas \cite{DBC08}. 
Semi arid here means 
that there is enough water to support {\em some} 
vegetation, but not enough water for a dense homogeneous vegetation. 
Such systems are often modeled in the form of two or more component 
reaction--diffusion systems 
for plant and water densities, with rainfall $R$ 
as the main bifurcation parameter, and the patterns 
are attributed 
to a positive feedback loop between plant density and water infiltration.   
Starting with a homogeneous equilibrium 
of high plant density for large $R$, stationary 
spatial patterns appear as $R$ is lowered, often following a 
universal sequence \cite{GRS14}.  This may lead to 
catastrophic, sometimes irreversible, regime shifts, where 
a patterned vegetation suddenly dies out completely, leaving 
a desert behind as $R$ drops below a certain 
threshold, or as some other parameter such as harvesting or 
grazing by herbivores 
is varied.  There is a rather large number 
of specific models, each agreeing with field observations in 
various parameter regimes, see e.g., \cite{meron01,riet01}, and 
\cite{riet04} for a review, 
or \cite{SBB09,meron13} and the references therein for further reviews 
including more recent work on 
early warning signals for desertification and other critical transitions, 
usually associated with subcritical bifurcations.   

A so far much less studied problem is the spatially distributed 
dynamic optimal control (OC) of vegetation systems 
by choosing harvesting or grazing by herbivores 
in such a space and time dependent way that some economic objective function 
is maximized. 
Following \cite{BX10} we consider an  infinite 
time horizon OC problem for a reaction--diffusion system for 
vegetation biomass and soil water, which is roughly based on 
\cite{riet01}. Related optimal control problems have also been considered 
in \cite{BX08,Xe10,Brocketal2013}, with the focus on the 
so called ``Optimal Diffusion Instability'' of flat 
canonical steady states (FCSS, see \S\ref{pssec} for OC related 
definitions), which 
similar to a Turing bifurcation yields the bifurcation of 
spatially {\em patterned} canonical steady states (PCSS) from FCSS. 
However, in these works, and 
in \cite{BX10}, only few PCSS have been actually calculated numerically, 
and no canonical paths, i.e., time dependent solutions of the canonical 
system. 
Finite time horizon cases recently have been considered in 
\cite{CP12, ACKT13, Apre14}, 
mostly focusing on theoretical aspects, and on 
problems with control constraints, which altogether gives a rather 
different setting than ours; see also Remark \ref{lrem}. 

Here we apply the numerical framework from \cite{GU15, U15p2} to the 
so called canonical system for the states and co--states, 
derived via \PMAXP. First we 
calculate (branches) of canonical steady states (CSS), including 
branches of PCSS, and in a second step 
canonical paths connecting to some CSS. Our main result is that 
in wide ranges of parameters, in particular for low rainfall $R$, 
FCSS and their canonical paths are not optimal, and that controlling 
the system to a PCSS yield a higher profit than uniform harvesting. 
{Thus, this seems to be the first example of a bio-economically 
motivated optimal control problem, where the global bifurcation structure 
of CSS has been computed in some detail, showing multiple CSS at 
fixed parameter values and with dominant non--spatially homogeneous 
steady states, and where moreover canonical paths to such PCSS have 
been computed, with significant gains in welfare.}  
We also compare these results to results for the same system 
with (initially) no external control, i.e., the system with private 
optimization, and find that the optimal control significantly 
increases the profit, and supports vegetation at significantly 
lower rainfall levels, {which again has important welfare implications. 
Remarkably, in our system the co-state 
of the vegetation can be identified with 
a tax for the private optimization problem, and thus, 
by solving the canonical system we find an optimal space and time 
dependent tax for the private optimization problem.} 

A standard reference on ecological economics or 
``Bioeconomics'' is \cite{Cl90}, including a very readable 
account, and applications, of \PMAXP{} in the context of 
ODE models; see also \cite{LW07}. 
A review of management rules for semi arid grazing systems including 
comparison to real data, and making plain the importance of such rules, 
is given in \cite{QB12}. However, the 
(family of) models discussed there consist of discrete time 
evolutions without spatial dependence, and with focus on the 
rainfall $R$ as a time--dependent stochastic parameter; here we 
consider a deterministic PDE model. Thus comparison between 
the two (classes of) models is difficult, but it should be interesting 
to include spatial dependence into the model in \cite{QB12}. 

Finally, in, e.g., \cite{Neu03}, \cite{DL09}, 
stationary 
spatial OC problems for a fishery model are considered, including 
numerical simulations, which correspond to our calculation of 
canonical steady states for our model \reff{oc1} below. 
The results of \cite{Neu03, DL09} show that for their models it is 
{\em economically optimal} to provide ``no--take'' marine reserves. 
This is similar to our finding of {\em optimal} patterned 
canonical steady states, but 
here we go beyond the steady case with the computation of optimal paths. 
This in particular tells us how to 
dynamically control the system to an (at least locally) 
optimal steady state. 
Moreover, even for the steady states there are important 
technical differences between the OC problems considered in 
\cite{Neu03, DL09} with {\em unique} positive canonical steady states, 
and the OC problems considered here, which  in relevant parameter regimes 
are distinguished by having {\em many} canonical steady states, of which several 
can be locally stable (see Remark \ref{defrem} for definitions, and, again, 
 Remark \ref{lrem}.)
  {Thus, management rules for our system are considerably 
more complicated than those in, e.g., \cite{Neu03,DL09}, as they 
must take the different CSS into account, and given an initial state, 
must decide to the domain of 
attraction of which CSS it belongs. Partly to illustrate this 
point, we also compute a Skiba (or indifference) point \cite{Skiba78} 
between different CSS in \S\ref{pskibasec}. }

In the remainder of this Introduction we explain the model and the use 
of \PMAXP{} to derive the canonical system (\S\ref{pssec}), explain 
the basic idea of the numerical method (\S\ref{sssec}), and summarize 
the results (\S\ref{srsec}). In \S\ref{rsec} we present the 
quantitative results, and in \S\ref{dsec} we give a 
brief further discussion. 
The method has been implemented in {\tt Matlab} as an ``add-on'' package 
{\tt p2pOC} to the continuation and bifurcation software {\tt pde2path} 
\cite{p2pure}, and 
the matlab functions and scripts  to run (most of) the simulations, the 
underlying libraries, manuals of the software, and some more demos   
can be downloaded at \cite{p2phome}. \\[3mm]
\bmip
{\bf Acknowledgment.} I thank D. Gra\ss, ORCOS Wien, 
for valuable comments on the economic terminology, and the anonymous 
referees for valuable questions and comments on the first version of 
the manuscript. 
\emip

\subsection{Problem setup}\label{pssec}
In dimensionless variables the model for vegetation and soil water 
with social optimization from \cite{BX10} is as follows. 
Let $\Om\subset\R^d$ be a one--dimensional (1D) or two--dimensional (2D) 
bounded domain, 
$v=v(x,t)$ the vegetation density at time $t\ge 0$ and space $x\in\Om$, 
$w=w(x,t)$ the soil water saturation, 
\bce
$E=E(x,t)$ the harvesting 
effort (control), $H(v,E)=v^\al E^{1-\al}$ the harvest, and
\ece 
\huga{\label{j0}
J(v_0,w_0,E)=\frac 1 {|\Om|}
\int_0^\infty\er^{-\rho t}\int_\Om J_c(v(x,t),E(x,t))\dd x\dd t
}
the (spatially averaged) objective function, where $\rho>0$ is the discount 
rate, and 
\huga{
J_c(v,E)=pH(v,E)-cE
} 
is the (local) current value profit. This profit depends on the price $p$, 
the costs $c$ for harvesting, and $v$, $E$ in a classical Cobb--Douglas form 
with elasticity parameter $0<\al<1$. 
 The problem reads 
\begin{subequations}\label{oc1}
\hual{
V(v_0,w_0)&=\max_{E(\cdot,\cdot)} J(v_0,w_0,E), \quad\text{where}\\
\pa_t v&=d_1\Delta v+[gwv^\eta-d(1+\del v)]v-H(v,E),\\
\pa_t w&=d_2\Delta w+R(\beta+\xi v)-(r_u v+r_w)w, \\
(v,w)|_{t=0}&=(v_0,w_0). 
}
\end{subequations} 
The parameters and 
default values are explained in Table \ref{tab1}, together with 
brief comments on the terms appearing in \reff{oc1}. One essential 
feature of (\ref{oc1}b,c) is the positive feedback loop between 
vegetation $v$ and water $w$. Clearly, the vegetation 
growth rate $gwv^\eta$ increases with $w$, but the vegetation 
$v$ also has a positive effect on water infiltration, 
for instance due to loosened soil, modeled by 
$R\xi v$ in (\ref{oc1}c).   
\brem\label{grazrem}{\rm In \cite{BX10}, $E$ is also refered to as 
the grazing by herbivores, and \reff{oc1} is called a semi arid 
grazing system. This seems somewhat oversimplified, because from a practical 
point of view, controlling in detail the movement and grazing behavior 
of, e.g., cattle, certainly is a more complicated problem than 
controlling genuine harvesting. Thus, henceforth we stick to calling $E$ 
the harvesting effort.}
\eex\erem 

We complement 
\reff{oc1} with homogeneous Neumann boundary conditions (BC) 
$\pa_\nu v=0$ and $\pa_\nu w=0$ on $\pa\Om$, where $\nu$ is 
the outer normal. The discounted 
time integral in \reff{j0} is typical for economic (here bioeconomic) 
problems, where ``profits now'' weight more than mid or far 
future profits. More specifically, $\rho$ corresponds to a long-term 
investment rate. We normalize $J_{ca}$ by $|\Om|$ for easier 
comparison between different domains and space dimensions. 
In (bio)economics, the control $E$, chosen externally by a ``social 
planner'' to maximize the social value $J$, is often called social 
control, as opposed to private optimization, see \reff{pde2} below. 
Finally, 
the $\max$ in (\ref{oc1}a) runs over all {\em admissible} controls $E$; 
essentially this means that $E\in L^\infty([0,\infty)\times\Om,\R)$, 
where moreover implicitly we have the control constraint $E\ge 0$, and 
state constraint $v,w\ge 0$ for the associated solutions of (\ref{oc1}b,c). 
However, in our simulations these constraints will always 
naturally be fulfilled, 
i.e., {\em inactive}, see also Remark \ref{lrem}. 
{\small 
\begin{table*}\bce\begin{tabular}{|p{10mm}|p{92mm}|p{38mm}|} 
\hline
param.&meaning&default values\\
\hline
$g,\eta$&coefficient and exponent in plant growth rate 
$gwv^\eta$&$g=0.001, \eta=0.5$\\
$d,\del$&coefficients in plant death rate 
$d(1+\del v)$&$d=0.03, \del=0.005$\\
$\beta,\xi$&coefficients in the infiltration function $\beta+\xi v$,
&
$\beta=0.9, \xi=0.001$\\
$R$&rainfall parameter, used as main bifurcation parameter&
between 4 and 100\\
$r_u, r_w$&water uptake and evaporation parameters in the water loss 
rate $r_u v+r_w$&$r_u=0.01, r_w=0.1$\\
$d_{1,2}$& diffusion constants for vegetation and water (resp.)& 
$d_1=0.05$, $d_2=10$\\\hline
$\rho$&discount rate&$\rho=0.03$\\
$c,p,\al$&(economic) param.~in the harvesting 
\mbox{$H(v,E)=v^\al E^{1-\al}$} and in the value $J_c(v,E){=}pH(v,E){-}cE$&
$c=1, p=1.1, \al=0.3$\\
&($\rho$ and $p$ used as bifurcation param.~in \S\ref{rhosec})&\\
\hline
\end{tabular}
\caption{Dimensionless parameters and default values in \reff{oc1}; 
see \cite{BX10} for further comments on the modeling. In particular, 
following \cite[\S4.2]{BX10} we have a rather larger $d_2$. \label{tab1}}
\ece 
\end{table*}
}

Introducing the costates $(\lam,\mu)=(\lam,\mu)(x,t)$ and the 
(local current value) Hamiltonian
\hual{
    \CH(v,w,\lam,\mu,E)=J_c(v,E)
&+\lam\bigl[d_1\Delta v+(gwv^\eta-d(1+\del v))v-H\bigr]\notag\\
&+\mu\bigl[d_2\Delta w+R(\beta+\xi v)-(r_u v+r_w)w\bigr], 
\label{hammax2}
}
by \PMAXP{} for $\tilde{\CH}=
\int_0^\infty \er^{-\rho t} \ov{\CH}(t)\dd t$ with the spatial integral 
\huga{\label{fullH} 
\ov{\CH}(t)=\int_\Om \CH(v(x,t),w(x,t),\lam(x,t),\mu(x,t),E(x,t))
\dd x, 
}
an optimal solution $(v,w,\lam,\mu)$ has to solve the canonical system (CS) 
\begin{subequations}\label{cs}
\hual{
\pa_t v&=\pa_\lam\CH=d_1\Delta v+[gwv^\eta-d(1+\del v)]v-H,\\
\pa_t w&=\pa_\mu\CH=d_2\Delta w+R(\beta+\xi v)-(r_u v+r_w)w,\\
\pa_t\lam&=\rho\lam-\pa_v\CH=\rho\lam-p\al v^{\al-1}E^{1-\al}-
\lam\bigl[g(\eta+1)wv^\eta-2d\del v-d-\al v^{\al-1}E^{1-\al}]\\\notag 
&\qquad\qquad\qquad-\mu(R\xi-r_uw)-d_1\Delta \lam, 
\\
\pa_t\mu&=\rho\mu-\pa_w\CH=\rho\mu-\lam g v^{\eta+1}+\mu(r_u v+r_w)-d_2\Delta \mu, 
}
where $E=\argmax_{\tilde{E}}\CH(v,w,\lam,\mu,\tilde{E})$, which is obtained from 
solving $\pa_E \CH=0$ for $E$, giving 
\huga{
E=\left(\frac{(p-\lam)(1-\al)}{c}\right)^{1/\al}v. 
}
The costates $(\lam,\mu)$ also fulfill zero flux BC, and 
derivatives like $\pa_v \CH$ etc are taken 
variationally, i.e., for $\ov{\CH}$. For instance, 
for $\Phi(v,\lam)=\lam\Delta v$ we have $\ov{\Phi}(v,\lam) 
=\int_\Om \lam\Delta v\dd x
=\int_\Om (\Delta\lam)v\dd x$ by Gau\ss' theorem, hence 
$\delta_v \ov{\Phi}(v,\lam)[h]=\int (\Delta\lam) h\dd x$, and 
by the Riesz representation theorem we 
identify $\delta_v \ov{\Phi}(v,\lam)$ and hence $\pa_v\Phi(v,\lam)$ 
with the multiplier $\Delta\lam$. Moreover, 
we used the so called intertemporal 
transversality conditions 
\huga{\label{tcond}
\lim_{t\ra\infty}\er^{-\rho t}\int_\Om v(x,t)\lam(x,t)\dd x=0 
\text{ and } \lim_{t\ra\infty}\er^{-\rho t}\int_\Om w(x,t)\mu(x,t)\dd x=0, 
}
which is justified since we are only interested in bounded solutions. 
\end{subequations}

For convenience setting 
\huga{
u:=(v,w,\lam,\mu):\Om\times [0,\infty)\ra\R^4, 
}
we collect (\ref{cs}a--e) and the boundary conditions into
\begin{subequations}\label{cs2}
\hual{\label{cs2a}
\pa_t u&=-G(u):=\CD\Delta u+f(u), \quad \CD=\bpm d_1&0&0&0\\0&d_2&0&0\\
 0&0&-d_1&0\\0&0&0&-d_2\epm,\\
\pa_\nu u&=0\text{ on } \pa\Om, \quad (v,w)|_{t=0}=(v_0,w_0). \label{cs2b}
} 
\end{subequations}
In \reff{cs2b} we only have initial conditions for the states, 
i.e., half the variables, and 
 \reff{cs2a} is ill-posed as an initial value problem 
 due to the backward diffusion in $(\lam,\mu)$. Thus, below 
we shall further restrict the transversality condition \reff{tcond} 
to requiring that $u(t)$ converges to a steady state, i.e.~a solution 
of 
\huga{\label{cs3}
G(u)=0,\quad  \pa_\nu u=0\text{ on } \pa\Om. 
}

\brem\label{defrem}{\rm {\bf (Definitions and notations)} 
A solution $u$ of the canonical system \reff{cs2} is called a 
\emph{canonical path}, 
and a solution $\uh$ of \reff{cs3} is called a 
\emph{canonical steady state (CSS)}. With a slight abuse of notation 
we also call $(v,w,E)$ with $E$ given by (\ref{cs}e) a canonical 
path, suppressing the associated co-states $\lam,\mu$. 
In particular, if $\hat u$ is a CSS, so is $(\hat v,\hat w,\hat E)$. 
A CSS $\hat u$ is called {\em flat} if it is spatially homogeneous, 
and {\em patterned} otherwise, and we use the acronyms 
FCSS and PCSS, respectively. 
For convenience, given $t\mapsto u(t)$ 
we also write 
\huga{\label{jredef}
J(u):=J(v_0,w_0,E), 
} 
with $(v_0,w_0)$ and $E$ (via (\ref{cs}e)) taken from $u$. 
A canonical path $u$ (or $(v,w,E)$) 
is called {\em optimal} if there is no canonical 
path starting at the state values $(v(0), w(0))$  
and yielding a higher $J$ than $J(u)$. As a special case, a 
CSS $\uh=(\hat v,\hat w, \hat \lam, 
\hat \mu)$ is called optimal if there is no canonical 
path starting at $(\hat v, \hat w)$  and yielding a higher $J$ 
than $J(\hat v,\hat w, \hat E)$. We use the acronyms OSS for any optimal 
CSS, and FOSS 
and POSS for optimal flat or patterned CSS $\uh$, respectively. 
An OSS $\uh$ 
is called {\em locally stable}   if 
for all admissible $(v_0,w_0)$ close to $(\hat v, \hat w)$ 
there is an optimal path $u$ with 
and $\lim_{t\ra\infty}u(t){=}\uh$.\footnote{For ``close to'' and $\lim$ 
we may use, e.g., the $H^1(\Om)$ norm, but since all solutions will be 
smooth, for instance as solutions of semilinear elliptic systems with 
smooth coefficients, we decided to omit the introduction of function 
spaces here. Similarly,  $(v_0,w_0)$ ``admissible'' should be read as 
$v_0(x),w_0(x)\ge 0$ for $x\in\Om$.}  
Similarly, $\uh$ is called {\em globally stable} if for all 
admissible $(v_0,w_0)$ the associated optimal path has 
$\lim_{t\ra\infty}u(t){=}\uh$. 
Finally, the 
{\em domain of attraction} of a locally (or globally) 
stable OSS $\uh$ is defined as the set of all $(v_0,w_0)$ such that 
the associated optimal path $u$ fulfills $\lim_{t\ra\infty}u(t){=}\uh$. 
See also \cite{GU15} for more formal definitions, and further 
comments on the notion of optimal system, and, e.g., 
the transversality condition \reff{tcond}. 
}\eex \erem

\brem\label{lrem}{\rm 
For background on OC in a PDE setting see for instance 
\cite{Tr10} and the references therein, or specifically \cite{RZ99, RZ99b, LM01} and 
\cite[Chapter5]{AAC11} for \PMAXP{} 
for OC problems for semilinear parabolic state evolutions. However, 
these works are in a 
finite time horizon setting, often the objective function 
is linear in the control, and there are state or control constraints. 
For instance, denoting the control by $k=k(x,t)$, often $k$ is chosen 
from some bounded interval $K$, and therefore is not obtained from 
the analogue of (\ref{cs}e), but rather takes values from $\pa K$, which 
is usually called bang-bang control.  
In, e.g., \cite{CP12, ACKT13}, some specific models of this type have been 
studied in a rather theoretical way, i.e., the focus is on 
deriving the canonical system and showing well-posedness and 
the existence of an optimal control. \cite{Apre14} additionally 
contains numerical simulations for a finite time horizon control--constrained 
OC problem for a three species spatial predator-prey system, again 
leading to bang--bang type controls.  See also \cite{NPS11} 
and the references therein for numerical methods for (finite time horizon)  
constrained parabolic optimal control problems. 

Similarly, the (stationary) fishery problems in \cite{Neu03, DL09} 
come with harvesting (i.e.~control) constraints. Moreover, in contrast 
to our zero--flux boundary conditions (\ref{cs2}b) 
Dirichlet boundary conditions are imposed. 
In \cite{KXL15} the models from \cite{Neu03, DL09} are extended to Robin boundary conditions, 
and a finite time horizon, with a discounted profit of the form 
$J=\int_0^T \int_\Om p\er^{-\rho t}h(x,t)u(x,t)\dd x\dd t$, 
where $p,\rho>0$ denote the price and discount rate, $h$ is the harvest, 
and $u$ the (fish) population density, which fulfills a rather general 
semilinear parabolic equation including advection. The first focus is 
again on well--posedness and the first order optimality conditions, and 
numerical simulations are presented for some specific model choices, 
illustrating the dynamic formation and evolution of marine reserves. 
However, the setting again is quite different from ours, due to the 
finite time horizon, and since $J$ is linear in $h$ and $u$, and 
since consequently there are constraints on $h$, leading to (unique) 
bang--bang controls. 

Here we do not consider (active) 
control or state constraints, and no terminal time, but the infinite 
time horizon. 
Our models and method are motivated by \cite{BX08, BX10}, which also 
discuss \PMAXP{} in this setting.  We do not aim at theoretical results, but 
rather consider 
\reff{cs2} after a spatial discretization as a (large) ODE problem, 
and essentially treat this using the notations and ideas from 
\cite{BPS01} and \cite[Chapter 7]{grassetal2008}, and the algorithms from 
\cite{GU15, U15p2}, to numerically compute optimal paths.}\eex
\erem 

\subsection{Solution steps}\label{sssec}
Using the canonical system \reff{cs2} 
we proceed in two steps: 
first we compute (branches of) CSS, and 
second we solve the 
``connecting orbit problem'' of computing canonical paths connecting 
to some CSS.  
Thus  we take a 
broader perspective than aiming at computing just one optimal control, 
given an initial condition $(v_0,w_0)$, which without further information 
is an ill-posed problem anyway. Instead, our method aims to give a 
somewhat global 
picture by identifying the (in general multiple) 
optimal CSS and their respective domains 
of attraction, as follows: 
\bci 
\item[(i)] To calculate 
CSS (which automatically fulfill \reff{tcond}) 
we set up \reff{cs3} as a bifurcation problem 
 and use the package \p2p\ \cite{p2pure, p2p2} 
to find a branch of FCSS, from which various branches of PCSS bifurcate. 
We focus on the rainfall parameter $R$ as the main bifurcation parameter, 
but in \S\ref{rhosec} 
also briefly discuss the dependence on the discount rate $\rho$ and 
the price $p$ as bifurcation parameters.
By calculating in parallel the (spatially averaged) current value profits 
\huga{\label{jca}
J_{c,a}(v,E):=\frac 1 {|\Om|} \int_\Om J_c(v(x),E(x))\dd x
}
of the CSS we can moreover immediately find which of the CSS maximize $J_{c,a}$ 
and hence $J=J_{c,a}/\rho$ amongst the CSS. 
\item[(ii)] In a second step we calculate canonical paths ending 
at a CSS (and 
often starting at the state values of a different CSS), and the objective 
values of the canonical paths. 
\eci 
Using a Finite Element Method (FEM) discretization, \reff{cs2} is 
converted into the ODE system (with a slight abuse of notation) 
\huga{\label{dcs} 
M\ddt u=-G(u),
} 
where $u\in\R^{4n}$ is a large vector containing the nodal values 
of $u=(v,w,\lam,\mu)$ at the $n$ spatial discretization points 
(typical values are $n=30$ to $n=100$ in 1D, and $n=1000$ and larger in 2D), 
and $M\in \R^{4n\times 4n}$ is the so called mass matrix, which is large but sparse. In \reff{dcs} and the following we mostly suppress 
the dependence of $G$ on the rainfall $R$ (or the other parameters).
 For (i) we thus need to solve the problem $G(u)=0$, which can be considered 
as an elliptic problem after changing the signs in the equations 
$G_{3,4}(u)=0$ for the costates. 

For (ii) we choose a suitable truncation time $T>0$ and 
replace the transversality condition \reff{tcond} by the 
condition 
\huga{ 
u(T)\in W_s(\hat u) \text{ and } \|u(T)-\hat{u}\| \text{ small}, 
}
where $W_s(\hat{u})$ is the stable manifold of the CSS 
$\hat{u}$ for the finite dimensional approximation 
\reff{dcs} of \reff{cs}. In practice we use 
\huga{
u(T)\in E_s(\hat u) \text{ and } \|u(T)-\hat{u}\| \text{ small},
}
where $E_s(\hat u)$ is the stable eigenspace of $\hat u$, i.e., 
the linear approximation of $W_s(\hat{u})$ at $\uh$. At $t=0$ we already 
have the  boundary conditions 
\huga{\label{ic1} 
(v,w)|_{t=0}=(v_0,w_0)\in\R^{2n} 
}
for the states. 
To have a well-defined two point boundary value problem in time 
we thus need 
\huga{\label{spp} \dim E_s(\hat u)=2n.}
Since the (generalized, in the sense of $M$ on the left hand side of \reff{dcs}) 
eigenvalues of the linearization $-\pa_u G(\hat u)$ of \reff{cs2} around 
$\uh$ are always symmetric around $\rho/2$, see \cite[Appendix A]{GU15}, 
we always have $\dim E_s(\hat u)\le 2n$. The number 
\huga{\label{ddef}
d(\hat u)=2n-\dim E_s(\hat u)
} 
is called the {\em defect} of $\hat u$, a CSS $\hat u$ with $d(\hat u)>0$ 
is called {\em defective}, and if $d(\hat u)=0$, then 
$\uh$ has the so called {\em saddle point property} (SPP). Clearly, these 
are the only CSS such that for all 
$(v_0,w_0)$ close to $(\hat v,\hat w)$ we may expect a solution 
for the connecting orbits problem 
\huga{\label{tbvpCS}
M\ddt u=-G(u), \quad (v,w)|_{t=0}=(v_0,w_0), \quad 
u(T)\in E_s(\hat u), \text{ and $\|u(T)-\uh\|$ small.}
} 
See \cite{GU15} for further comments on the significance 
of the SPP \reff{spp} on the discrete level, and its (mesh-independent) 
meaning for the canonical system as a PDE. 

\subsection{Results}\label{srsec}
The bifurcation diagram (i) for \reff{cs3} turns out to be quite rich, 
already over small 1D domains. Thus we mostly focus on 1D, 
but we include a short 2D discussion in \S\ref{2dsec}. 
Details will be given in \S\ref{rsec}, 
but in a nutshell we have: In pertinent $R$ regimes there are 
many CSS, but most of them do not fulfill the SPP, and most of 
these that do fulfill the SPP are not optimal. On the other hand, 
in particular at low $R$ there are locally stable POSS 
(patterned optimal steady states). 
Before further commenting 
on this, we briefly review results for the so called uncontrolled case.  

In \cite[\S2]{BX10} it is shown that the case of private 
objectives, where economic agents (ranchers) are located at each site $x$, 
and each one maximizes his or her private profits, leads to the 
system 
\begin{subequations}\label{pde2}
\hual{
\pa_t v&=d_1\Delta v+(gwv^\eta-d(1+\del v)-A)v,\\
\pa_t w&=d_2\Delta w+R(\beta+\xi v)-(r_u v+r_w)w, 
}
\end{subequations}
i.e., the harvest is $H=Av$, $A>0$ a fixed parameter. In detail, 
 ranchers with certain property rights individually 
maximizing $\pi(v,E)=p v^\al E^{1-\al}-cE$  leads to 
\huga{\label{poE} 
\text{$E=\ga v$ with 
$\ga=\ds\left(\frac{p(1-\al)}{c}\right)^{1/\al}$ and hence $A=
\ga^{1-\al}$.}
}
The same happens in the so called open access case, where agents may 
harvest freely, giving 
$E=\hat\ga v$ with $\hat\ga=\ds \left(\frac{c}{p}\right)^{-1/\al}$ and 
hence $A=\hat\ga^{1-\al}$. On the other hand, \reff{pde2} can also be seen 
as a ``undisturbed Nature'' case with modified vegetation death rates 
$\tilde{d}=d+A$ and $\tilde{\del}=\frac d {d+A}$. 
For the economic parameters $(c,p,\al)=(1,1.1,0.3)$ from Table \ref{tab1} 
we obtain $A=0.543$, which is rather large compared to the original 
$d=0.03$ from Table \ref{tab1}. 

The bifurcation picture for steady states of \reff{pde2} is 
roughly similar to that of \reff{cs3}, but there are also significant 
differences, and we altogether summarize our results as follows:  
\bci
\item[{\bf (a)}] For large $R$ there is a ``high vegetation''  FCSS, 
which is a globally stable FOSS (flat optimal steady state). 
\item[{\bf (b)}] For smaller $R$ the FCSS from {\bf (a)} loses optimality, 
and there bifurcate branches of (locally stable) POSS 
(Patterned Optimal Steady States). 
In particular, we can calculate canonical paths from 
the state values of the FCSS to some POSS which increase 
the profit (up to 40\%, see \ref{p0val}). 
\item[{\bf (c)}] The uncontrolled flat steady states (FSS) of \reff{pde2} 
only exist 
for much larger $R$ values than the FCSS. 
At equal $R$, the profit $J$ (or equivalently the 
discounted value $J/\rho$) of the FSS is much lower (e.g., one tenth, see 
Table \ref{tab3}, bottom center)  
than the value of the FCSS of \reff{cs3}. 
\item[{\bf (d)}] For the initial value problem \reff{pde2} we may consider the 
stability of steady states, while CSS at best have the SPP. 
It turns out that the FSS branch loses stability at a much larger 
value of $R$ to a patterned steady state of \reff{pde2},  
than the FOSS loses optimality in 
the optimal control problem \reff{oc1}. Thus, in the uncontrolled problem 
pattern formation sets in at larger $R$. 
\eci
Roughly speaking, {\bf (b)} means that at low $R$ it is advantageous 
to restrict harvesting to certain areas. This is similar 
to the marine reserves in the fishery models in \cite{Neu03, DL09}, 
and not 
entirely surprising, as it is well known both in models 
and in field studies of semi arid systems 
that also in uncontrolled systems low rainfall 
levels lead to 
patterned (or patchy) vegetation, which often can be sustained at 
lower rainfall levels than a uniform state \cite{meron01,riet04,meron13};  
this is precisely what happens here as well. 
However, the quantitative differences between steady states of 
\reff{cs2} and \reff{pde2} are quite significant, in particular in 
the sense {\bf (c),(d)}: The controlled system sustains vegetation 
at much lower $R$ values, and at equal $R$ it yields a much 
higher profit than \reff{pde2}.  Moreover, the computation of 
canonical paths from some initial state (often a CSS) to some OSS 
yields precise information  {\em how} to go to the OSS to maximize 
the objective function. 

{\brem\label{taxrem}{\rm In economics, 
the co-states $\lam$ and $\mu$ are also called {\em shadow prices}, 
which are sometimes 
difficult to interpret; here, however, the shadow price $\lam$ has a 
nice interpretation as an optimal tax for private optimization, as follows. 
Introducing a tax $\tau$ per unit harvest, i.e., 
setting 
\def\pit{\tilde{\pi}}
\huga{\label{potax1}
\pit(v,E)=(p-\tau)v^\al E^{1-\al}-cE,
}
we obtain 
\huga{\label{poEt} 
E=\left(\frac{(p-\tau)(1-\al)}{c}\right)^{1/\al}v
}
instead of \reff{poE}. Thus, after solving the optimization 
problem \reff{oc1} via \PMAXP, 
and thus in particular computing the shadow prize 
$\lam$ of the vegetation constraint along an optimal path, comparing \reff{poEt} 
to (\ref{cs}e) we see that 
private optimization maximizes $J$ from \reff{j0}, if the tax $\tau$ 
is set in an in general time and space dependent way as $\tau=\lam$. 
}
\eex\erem }

As already said, our method follows \cite{GU15}, where we study optimal 
controls for a model of a shallow lake ecology/economy, given by 
a scalar parabolic PDE. However, the results are rather different. 
First of all, without control, i.e., for a fixed control 
equal to some parameter, the scalar PDE in \cite{GU15} shows no 
pattern formation: the patterns in \cite{GU15} are only due 
to the control, which may be called ``control induced pattern formation'' 
or, as in \cite{BX08, BX10}, ``optimal diffusion instability'' (ODI). 
However, the parameters in \cite{GU15} have to 
be carefully fine--tuned to obtain POSS, which moreover are only locally stable, 
 see also \cite{grass2014}. 
Here, to obtain POSS we need no fine--tuning of parameters. 

From the methodical and algorithmic point of view, our results for \reff{oc1} 
illustrate that our two--step approach is well suited to deal with 
non--uniqueness of CSS in nonlinear PDE optimal control problems, 
and the typically associated multiple canonical paths and multiple local maxima 
of the objective function. See also \S\ref{dsec} 
for further discussion of efficiency. 
\section{Numerical simulations}\label{rsec}
\subsection{1D canonical steady states bifurcating from 
the FCSS branch} 
\label{1dcss-sec}
The bifurcation scenario 
for the stationary problem $G(u)=0$ can be studied conveniently 
with \p2p. First we concentrate on the 1D case $u=u(x)$, $x\in(-L,L)$, 
where the domain length must be chosen in such a way to capture 
pertinent instabilities of the FCSS branch. In \cite{BX08, BX10} conditions 
for pattern forming instabilities  
in terms of the Hamiltonian $\CH$ and its derivatives at a 
FCSS are given. These are similar to the well known 
Turing space conditions \cite{Murray89book}, and moreover allow 
the calculation of the critical wave--number $k_c$ of 
the bifurcation patterns. For instance, at $R=5$ \cite{BX10} 
find $(k_-,k_+)\approx (0.146, 1.455)$ for the band of 
unstable wave numbers at the FCSS. 

If one is interested in accurately capturing  the first bifurcation, 
then one should 
either fit the domain to the (wave number of the) first instability (see, e.g., 
\cite{uwsnak14} for examples), or use a very 
large domain, which gives a rather dense set of allowed wave numbers. 
However, for simplicity, and with the (expensive) $t$--dependent 
canonical paths in mind, here we do not want 
to use a very large domain, and, moreover, rather take the point 
of view that the domain comes with the model. Thus, we 
do not want to be too 
precise on fitting the domain to the first instability over an 
infinite domain, and simply choose $L=5$. Of course, increasing the domain size 
(certainly in integer multiples of $L$) will only increase the 
number of patterns and bifurcations, and on the other hand 
there is a critical minimal 
domain size below which no patterns exist.

\begin{figure}[ht]
\hs{-15mm}\ig[width=1.2\textwidth, height=130mm]{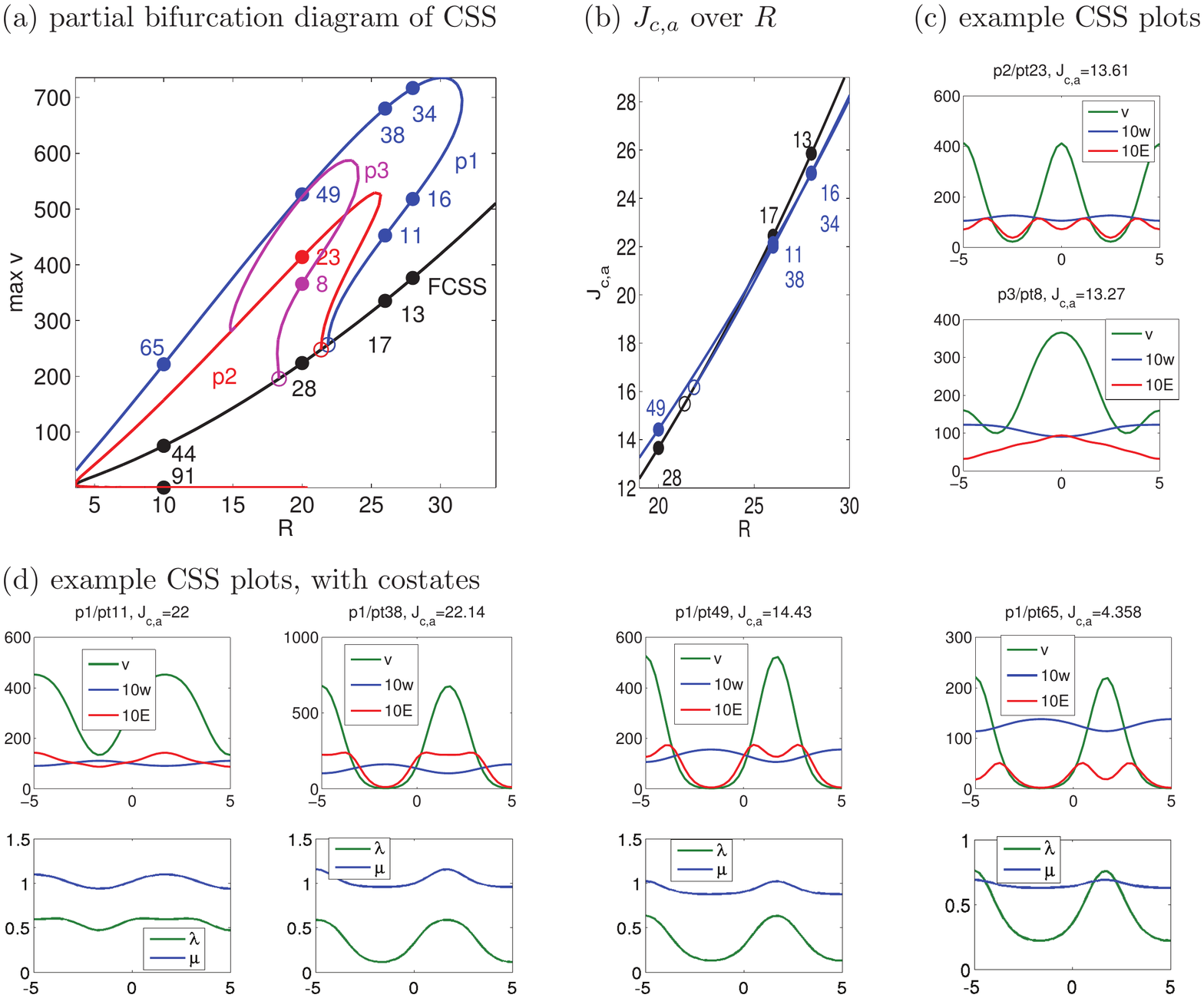}
\vs{-8mm}
\caption{{\small (a) (partial) bifurcation diagram of CSS in 1D, including 
a (small) selection of bifurcation points indicated by $\circ$. The 
labeled points are tabulated in Table \ref{tab3}, and a selection 
of these solutions together with $E$ (and the co--states for the 
p1--branch) are plotted 
in (c),(d). 
(b) shows the two branches 
FCSS and p1 with colors as in (a) in a $J_{c,a}$ over $R$ diagram, 
which allows to identify which CSS maximizes $J_{c,a}$ amongst the CSS 
(all other branches from (a) have a significantly lower $J_{c,a}$). 
\label{f1}}}
\end{figure}

In order to present our results in a domain independent way 
we give averaged quantities such as $J_{c,a}$, see \reff{jca}, and 
\huga{\label{sprdef} 
\spr{v}:=\frac 1 {|\Om|} \int_\Om v(x)\dd x, \quad 
\text{ and so on}. 
}
Figure \ref{f1}(a) shows a basic bifurcation diagram of 1D CSS. We start 
with a FCSS at $R=34$ which can be easily found from an initial guess 
$(v,w,\lam,\mu)=(400, 10, 0.5, 1)$ followed by a Newton 
loop.\footnote{(\ref{cs}a,b) also has 
the trivial solution branch $(v,w)=(0,R\beta/r_w)$, which yields the 
trivial branch FCSS$_0$ with 
$(v,w,\lam,\mu)=(0,R\beta/r_w,0,0) \text{ (and hence $E=0$),}$ 
which however is of little interest here.}
See Table \ref{tab3} 
for numerical values of some of the FCSS found this way, and of other 
selected points in the bifurcation diagrams. Following the FCSS branch 
with decreasing $R$ we find a 
number of branch points to PCSS, and near $R=4$ 
we find a fold for the FCSS branch. The lower FCSS branch continues back 
to large $R$, but is not interesting from a modeling point 
of view. The upper FCSS branch has the SPP until the first bifurcation 
at $R=R_c\approx 21.5$, where a PCSS branch p1 with period $20/3$ 
bifurcates subcritically; see the example plots of solutions on p1. 
This is a pitchfork 
bifurcation, but here and in general we only follow the branch in one direction; 
the other direction is often related to the first by symmetry, e.g., spatial 
translation by half a period. 

\begin{table}[ht]\begin{footnotesize}
\bce
\begin{tabular}{|c|ccC{5mm}C{1mm}C{10mm}ccC{9mm}C{1mm}cc|}
\hline
point&FCSS/13&p1/16&p1/34&&FCSS/17&p1/11&p1/38&\mbox{\hs{-2mm}}FCSS/28&&p1/49&p2/23\\
R&28&28&28&&26&26&26&20&&20&20\\
$\spr{v}$&376.32&337.92&283.31&&335.36&311.63&252.46&223.59&&175.49&185.08\\
$\spr{w}$&9.25&10.37&13.17&&9.3&10.05&13.34&9.62&&13.28&11.64\\
$\spr{\lambda}$&0.59&0.53&0.35&&0.59&0.56&0.34&0.62&&0.34&0.45\\
$\spr{\mu}$&1.09&1.06&1.06&&1.04&1.02&1.03&0.87&&0.92&0.89\\
$J_{c}$&25.85&25.06&25.02&&22.45&22&22.14&13.66&&14.43&13.61\\
$d$&0&2&0&&0&2&0&2&&0&1\\
\hline 
point&FCSS/44&FCSS/91&p1/65&\vline&\mbox{\hs{-2mm}}FSS$^*$/121&p1$^*$/97&p6$^*$/30&\mbox{\hs{-2mm}}FCSS/105&\vline&q1/76&q6/69\\
R&10&10&10&\vline& 60&60&60&60&\vline&20&20\\
$\spr{v}$&75.08&0.21&71.33&\vline& 79.73&163.42&86.5&1304.5&\vline&183.14&151.22\\
$\spr{w}$&11.46&88.21&12.71&\vline& 65.51&39.21&61.4&10.06&\vline&12.02&14.54\\
$\spr{\lambda}$&0.68&0.9&0.42&\vline&NA&NA&NA&0.49&\vline&0.35&0.31 \\ 
$\spr{\mu}$&0.5&0.0006&0.65&\vline&NA&NA&NA&1.75&\vline&0.91 &0.91\\
$J_{c}$&3.51&0.002&4.36&\vline& 14.3&29.31&15.51&120.8&\vline&13.93&13.93\\
$d$&4&4&0&\vline& NA(u)&NA(s)&NA(u)&0&\vline&0&0\\
\hline
\end{tabular}
\caption{{\small Characteristics of selected points marked in 
Fig.~\ref{f1}(1D, top half and bottom left in the table), 
Fig.~\ref{fnc} (bottom center in the table, where $^*$ denotes values 
from the case of private optimization) and 
Fig.~\ref{f3} (2D case, bottom right in the table). NA for the ${}^*$ values 
means that these values are not defined; for the defect the additional 
u,s are used to indicate unstable vs stable solutions.}
\label{tab3}}
\ece \end{footnotesize}
\end{table}

The p1 branch has a fold at $R_f\approx 31$, after which it has the SPP down 
to a secondary bifurcation at $R_2\approx 9.4$. 
Near $R=3.1$ the p1 branch also 
has a second 
fold, after which it continues to larger $R$ as a branch of small amplitude  
patterns. (a) also shows the PCSS branches bifurcating from the second and 
third bifurcation points on the FCSS; interestingly, p3 connects back to 
p2 in a secondary bifurcation. However, except for the FCSS branch for 
$R>R_c$ and the p1 branch between the first fold and the first 
secondary bifurcation, no other branch has solutions with the SPP, 
cf.~\reff{spp}.

Figure \ref{f1}(b) shows the FCSS and p1 branches in a $J_{c,a}$ over 
$R$ bifurcation diagram. This already shows that at, e.g., $R=20$ 
(in fact, for $R$ smaller than about 24.4),  
solutions on p1 yield a larger $J_{c,a}$ than the FCSS, 
and thus appear as a candidates 
for POSS. From the applied point of view, probably the most 
interesting aspect of the solution plots in (c) and (d) is 
that after the first fold on the p1 branch 
the effort $E$ has local minima, not maxima, at the maxima of 
$v$. Instead, $E$ has maxima on the slopes near the maximum of $v$. 
{Taking into account (\ref{cs}e) and the co-state plots in 
the second row of (d), this can be attributed to the distinctive 
peaks in $\lam$ at the maxima of $v$. These peaks in the shadow price 
of the vegetation evolution illustrate that it is not optimal 
to harvest at the peaks in $v$ as this will strongly decrease future 
income. Also note that the (average) shadow prices $\spr{\lam}$ on 
the p1 branch after the fold are lower than on the FCSS 
branch at the same $R$, while at least at low $R$, $\spr{\mu}$ is higher 
on p1 than on FCSS. }

The vegetation patterns (p1 branch) 
survive for lower $R$ (up to $R_{\text{crit}}\approx 3.1$) than 
the FCSS branch ($R_{\text{crit}}\approx 3.7$). 
However, the difference is not large, and this bottom end of $R$ 
will not 
be our interest here, despite its significance for critical transitions. 
Instead, we are interested in the optimality 
of CSS for intermediate $R_{\min}<R<R_f$, with $R_{\min}=5$, say.  

\brem\label{unirem}{\rm Although our picture of CSS obtained above 
is already somewhat complicated, naturally it 
is far from complete. Firstly, we only followed the first three 
bifurcations from the FCSS branch, and secondly, there are (plenty of) 
secondary bifurcations on the branches p1, p2 and p3, which here 
we do not follow. 
In particular, given that the 1.5--modal (in, e.g., $v$) 
branch 
p1 maximizes $J$ amongst the CSS, a natural question is whether 
there also exist unimodal or 0.5--modal branches, which might given 
even higher $J$. The answer is (partly) yes: while we could not 
find a 0.5--modal branch, there is a unimodal 
branch p0, which bifurcates from p2 in a secondary bifurcation, 
or, more precisely, connects p1 and p2. See \S\ref{unisec} for details, 
where inter alia we study the bifurcation behavior in the price $p$. 
Moreover, p0 then maximizes $J$ amongst the CSS, and, loosely speaking, 
turns out to be a global maximum for \reff{oc1}. 

Nevertheless, until \S\ref{unisec}, 
for the sake of clarity we restrict to the primary 
branches which bifurcate from the FCSS when varying $R$. However, 
one should keep in mind that whatever method one uses to 
study optimization problems like \reff{oc1} it is always possible 
to be stuck with some local maxima, and to miss some global maximum. 
}
\eex 
\erem 

\subsection{Comparison to private optimization} 
\label{ncsec}
\begin{figure}
\bce
 \ig[width=0.99\textwidth, height=100mm]{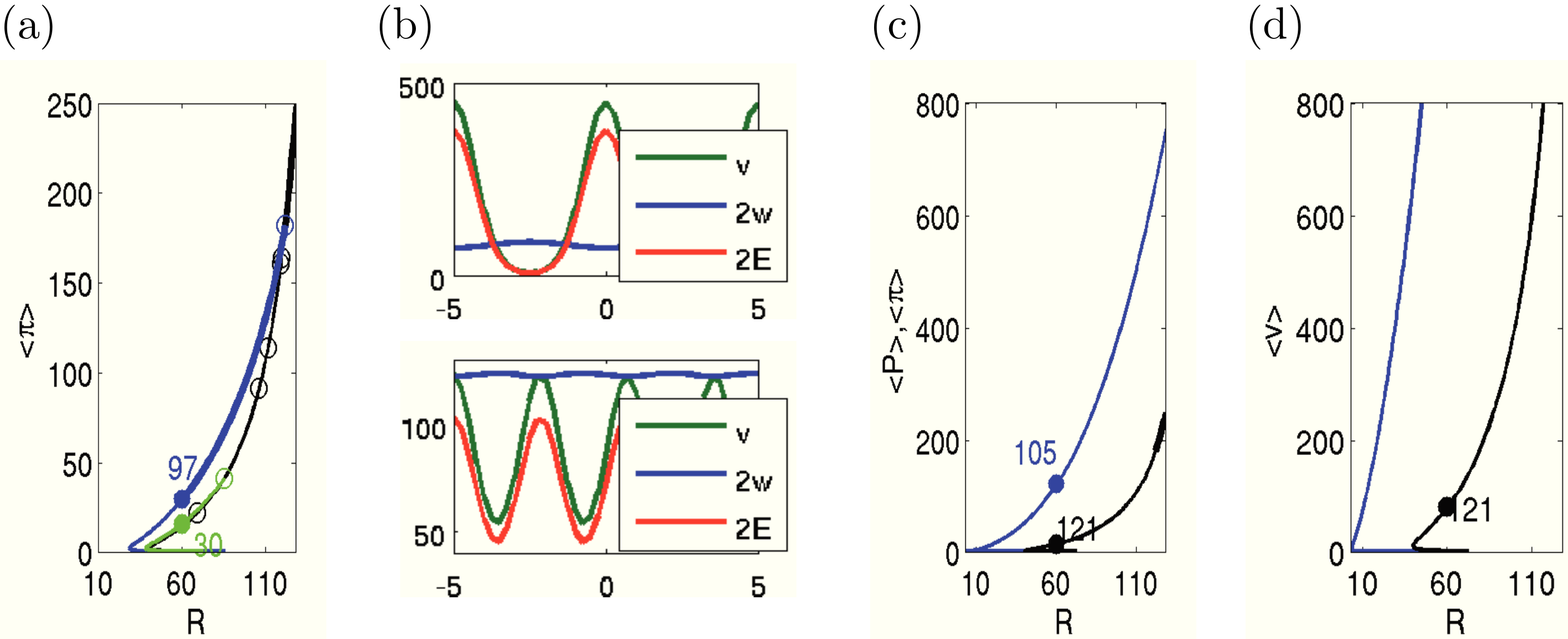}
\ece
\vs{-30mm}
\caption{{\small (a) Bifurcation diagram for the case of private optimization, 
where $\spr{\pi}$ denotes the average private profit (given by the same formula 
as the profit $P$ under control); 
the blue branch is the primary bifurcation {\tt p1nc}, the green is {\tt p6nc}. 
(b) example solutions from (a), p1nc/pt97 (top) and p1nc/pt30 (bottom). 
(c,d) comparison of the flat states without control (black) 
with the OC case (blue). \label{fnc}}}
\end{figure}

As already explained in the Introduction, private objectives, i.e., 
individual ranchers maximizing  $\pi(v,E)=p v^\al E^{1-\al}-cE$, 
leads to the system 
\begin{subequations}\label{pde12}
\hual{
\pa_t v&=d_1\Delta v+(gwv^\eta-d(1+\del v)-A)v,\\
\pa_t w&=d_2\Delta w+R(\beta+\xi v)-(r_u v+r_w)w, 
}
\end{subequations}
with $A=\ga^{1-\al}=0.543$ for the economic parameters 
$(c,p,\al)=(1,1.1,0.3)$ from Table \ref{tab3}. In this 
section we compare the bifurcation diagram in $R$ for steady states of 
\reff{pde12}, see Fig.~\ref{fnc}, to that for \reff{cs2} in Fig.\ref{f1}.

Roughly speaking, both are similar, but for \reff{pde12} the bifurcations to 
patterned steady states occur at larger $R$, 
and of course also have to be interpreted differently. First of 
all, we start the bifurcation diagram at $R=130$ with a dynamically 
stable flat steady state (FSS) of \reff{pde12}, which loses stability at $R_c\approx 122$ due to a 
supercritical pitchfork bifurcation of a branch p1nc of patterns with period 5. 
There are a number of further bifurcations from the FSS branch; 
as an example we give p6nc. The p1nc solutions lose stability in a secondary 
bifurcation near $R=61$ (not followed here), and eventually all branches 
undergo a fold between $R=36$ and $R=24$, and turn into small 
amplitude branches. 

Similar to the CSS case, here we also have $\spr{\pi({\rm p1nc})}>
\spr{\pi({\rm FSS})}$ 
i.e., the patterned states yield a higher (average) profit than the 
flat states
\footnote{although $J_{c,a}$ and $\spr{\pi}$ have different interpretations, 
in Table \ref{tab3} we use $J_{c,a}$ also for $\spr{\pi}$, as both are 
actually defined by the same expression}. 
In (c,d) we compare the FSS branch 
with the FCSS branch.   Besides again showing that the fold of the FCSS is at 
much lower $R$, and hence the socially 
controlled system supports a uniform vegetation 
down to much lower $R$, this also illustrates that, at given $R$, 
$\spr{P}$ and $\spr{v}$ are 
significantly higher on the FCSS branch.  Finally, 
\reff{pde12} has the trivial branch $(v,w)=(0,R\beta/r_w)$, 
which however again is of no interest to us. 

\subsection{Canonical paths}\label{1dcpsec}
\subsubsection{Main results}
Having computed CSS branches is only half of the program outlined above; 
we also need to solve the time dependent problem \reff{tbvpCS} to 
\bci
\item
given a point $(v_0,w_0)$ and a CSS $\hat{u}$, 
determine if 
there exists a canonical path from $(v_0,w_0)$ to $\hat u$; 
\item
compare  canonical paths $t\mapsto u(t)$ (or 
$(v,w,E)(t)$) 
to different $\hat u$, i.e.: compute and compare their economic values 
$J(u):=J(v_0,w_0,E)$, cf.~\ref{jredef}, 
to find {\em optimal} paths. 
\eci 
Assuming that the spatial mesh consists of $n$ points, we summarize the 
algorithm for the first point as follows, with more details 
given in \cite{GU15,U15p2}. First we compute $\Psi\in \R^{2n\times 4n}$ 
corresponding to the unstable eigenspace of $\hat u$ such that the 
right BC $u(T)\in E_s(\uh)$ is equivalent to $\Psi (u(T)-\uh)=0$. Then, to 
solve \reff{tbvpCS} we use a modification {\tt mtom} of the two-point BVP 
solver TOM \cite{MT04}, 
which allows to handle systems of the form $M\ddt u=-G(u)$, 
where $M\in\R^{4n\times 4n}$ is the mass matrix arising in the FEM 
discretization of \reff{cs}. A crucial step in solving 
(nonlinear) BVPs is  a good initial guess for a 
solution $t\mapsto u(t)$, and we combine {\tt mtom} with 
a continuation algorithm in the initial states, again see 
\cite[\S2.1]{GU15} 
for further discussion, and \cite{U15p2} for implementation details. 

For the second point we note that for a CSS $\uh$ we simply have  
$J(\uh)=J_{c,a}(\uh)/\rho$. 
Given a canonical path $u(t)$ that 
converges to a CSS $\uh$, and a final time $T$, we may then approximate 
\hual{
J(v_0,w_0,E)=&\int_0^T\er^{-\rho t}J_{c,a}(v(t,\cdot),E(t,\cdot))\dd t
+\frac {\er^{-\rho T}}{\rho}J_{c,a}(\uh). 
}

\begin{figure}[ht]\small
\hs{-20mm}
\ig[width=1.25\textwidth, height=140mm]{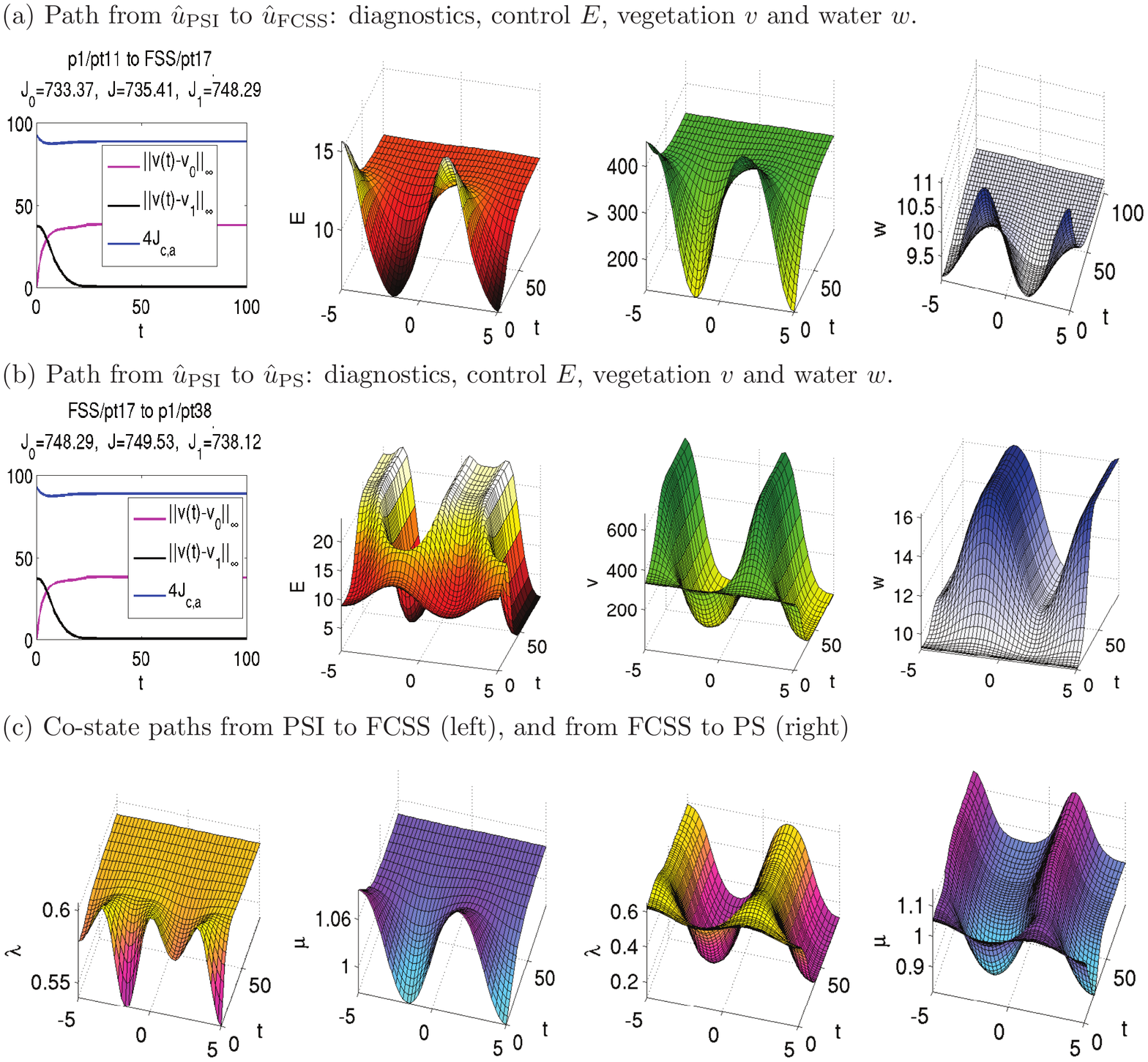}
\vs{-10mm}
\caption{{\small $R=26$. (a) convergence behavior, the 
current value profit, and obtained objective values for 
the canonical paths from p1/pt11 
to the FCSS (left), and $E,v,w$ on the path; (b) the same for the 
canonical path from the FCSS to p1/pt38. (c) co-state paths from (a),(b). }
\label{vf2}}
\end{figure}

In principle, given $\uh=(\hat v,\hat w, \hat \lam, \hat \mu)$ with 
$d(\uh)=0$ we could choose any 
$(v_0,w_0)$ in a neighborhood of $(\hat v,\hat w)$ (or globally, 
if $\uh$ is a globally stable OSS) and aim to find a canonical paths from 
 $(v_0,w_0)$ to $\uh$. However, the philosophy of our simulations 
rather is to start at the state values of 
some CSS, and see if we can find canonical paths to 
some other CSS which give a higher $J$. 
We discuss such canonical paths in decreasing $R$, starting with $R=26$ in 
Fig.~\ref{vf2}, and postponing the situation at $R=28$ to \S\ref{pskibasec}. 

\brem\label{prerem}{\rm a) Note again, that our discussion 
is based on the primary bifurcations in $R$ from the FCSS in Fig.~\ref{f1}, 
which misses a branch of unimodal CSS, cf.~Remark \ref{unirem} 
and \S\ref{unisec}.\\
b) Although only the state values $(v_0,w_0)$ are fixed as initial 
conditions for canonical paths as connecting orbits, in order not to 
clutter notations and language we write, e.g., $u_{\text{FCSS$\to$PS}}$ 
for a connecting orbit starting at the state values of 
$\uh_{\text{FCSS}}$ and connecting to $\uh_\text{PS}$. 
}\eex 
\erem

\paragraph{$R=26$.}
At $R=26$ we have two 
CSS with $d(\uh)=0$, namely $\uh_{\text{FCSS}}$ given by FSS/pt17, 
and $\uh_{\text{PS}}$ 
given by p1/pt38. Figure \ref{vf2} shows two canonical paths 
to these CSS. (a) shows $\|v(t)-v_0\|_\infty$ and $\|v(t)-v_1\|_\infty$ 
for the path from the ``intermediate'' patterned CSS $\uh_{\text{PSI}}$ 
given by p1/pt11 to 
$\uh_{\text{FCSS}}$. This indicates 
that and how fast the canonical path leaves the initial $(v_0,w_0)$ and 
approaches the goal $(v_1,w_1)=(\hat v,\hat w)$ (the differences 
in the second component $\|w(t)-w_*\|_\infty$ are always smaller). 
Moreover we plot $4J_{c,a}(t)$, illustrating 
that $J_{c,a}(t)$ does not vary much along the canonical path. 
However, the differences may accumulate. 

The values of the solutions are as follows. 
We readily have 
\huga{\label{cpv1} 
J(\uh_{\text{PSI}})=733.37<J(\uh_{\text{PS}})=738.12<J(\uh_{\text{FCSS}})=748.29, 
}
and for the paths $u_\text{PSI$\to$FCSS}$ (a,b), 
$u_\text{FCSS$\to$PS}$ (c), and 
$u_\text{PSI$\to$PS}$ (not shown) we obtain 
 \huga{\label{cpv2} 
J(u_\text{PSI$\to$FCSS})=735.51, \quad J(u_\text{PSI$\to$PS})=744.28, 
 \text{ and } J(u_\text{FCSS$\to$PS})=749.53. 
}
The result for $J(u_\text{PSI$\to$FCSS})$ 
seems natural, as controlling the system to 
a CSS with a higher value should increase $J$.  
However, the results for $J(u_\text{PSI$\to$PS})$ and 
$J(u_\text{FCSS$\to$PS})$ may at first seem counter--intuitive. 
In  $u_\text{FCSS$\to$PS}$ we go to a CSS with a smaller value, but 
the transients yield a higher profit for the path. 
In particular, this shows that 
{\em a CSS which maximizes $J$ amongst all CSS is not necessarily optimal}.
Similarly,  although $\uh_{\text{PS}}$ as a CSS has a lower value 
than $\uh_{\text{FCSS}}$, 
starting at $\uh_{\text{PSI}}$ it is advantageous to go to $\uh_{\text{PS}}$ 
rather than to $\uh_{\text{FCSS}}$. 

Due to folds in the continuation in the initial states, 
again see \cite{GU15} for details, we could not compute a path 
from $\uh_{\text{PS}}$ to $\uh_{\text{FCSS}}$. Thus we conclude 
that such paths do not exist, and (tentatively, see Remark \ref{prerem}) 
classify $\uh_{\text{PS}}$  as an at least locally stable 
POSS, with $\uh_{\text{FCSS}}$ and $\uh_{\text{PSI}}$ 
in its domain of attraction. 

The control to go from $\uh_{\text{PSI}}$ to $\uh_{\text{FCSS}}$ in (b) is  
intuitively clear: Increase/decrease $E$ near the maxima/minima of $v_0$.  Going 
from $\uh_{\text{FCSS}}$ to $\uh_{\text{PS}}$ in (c) 
warrants a bit more discussion: 
For a short transient, $E$ is reduced 
around the locations $x_2= -5$ and $x_4=5/3$ of the maxima of 
$\hat v_{\text{PS}}$. 
This is enough to give an increase of $v$ around $x_{2,4}$. 
However, under the given conditions this does not decrease soil water 
near $x_{2,4}$, i.e., the increased infiltration at larger $v$ 
dominates the higher uptake by plants. After this transient, 
$E$ increases near $x_{2,4}$, thus producing the higher $J$; 
 see also the discussion of Fig.~\ref{f2b} below. 
{As the behavior of $E$ follows from (\ref{cs}e), i.e., 
$E=\left(\frac{(p-\lam)(1-\al)}
{c}\right)^{1/\al}v$ for illustration 
we also plot $\lam,\mu$ for the paths in (a), (b) in Fig.~\ref{vf2}(c). }

\begin{figure}[ht]
\bce 

\vs{-25mm}
\ig[width=0.99\textwidth, height=105mm]{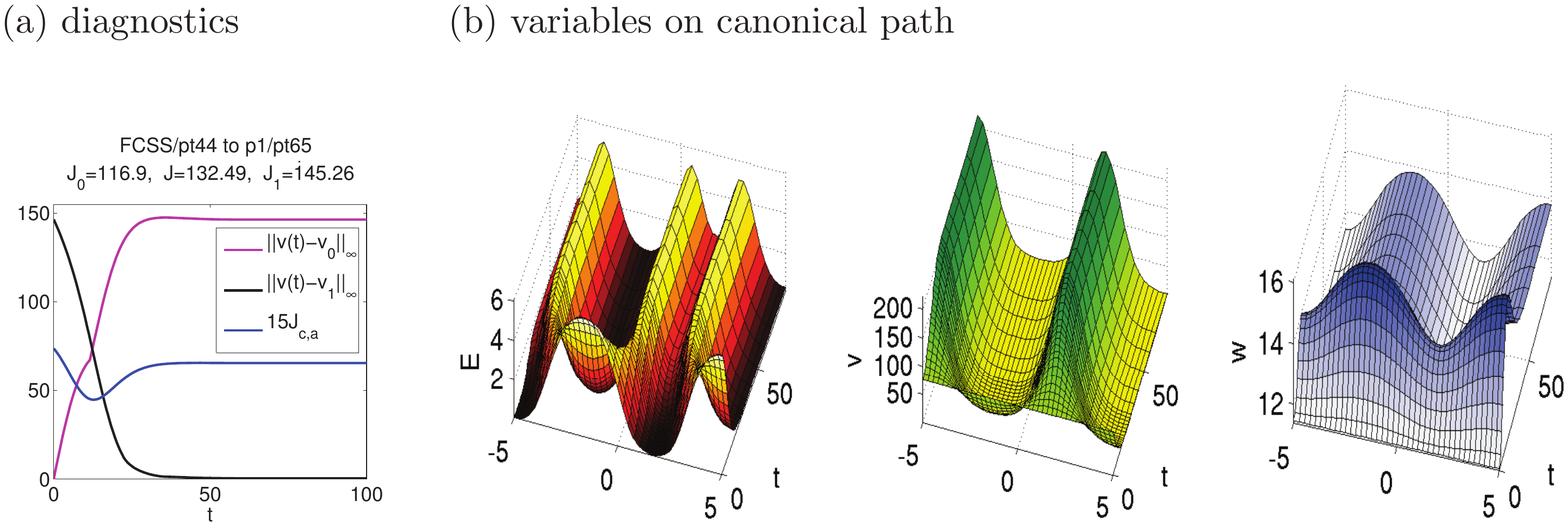}
\ece 
\vs{-35mm}
\caption{{\small The canonical path from FCSS to PS at $R=10$. 
\label{f2b}}}
\end{figure}

\paragraph{Smaller $R$.}
For  $R{<}R_c{\approx}21.5$, in Fig.\ref{f1} the only CSS with the SPP is the 
patterned state $\uh_{\text{PS}}(R)$. In Fig.~\ref{f2b} we 
focus on the case $R{=}10$, and here only remark that the results are 
qualitatively similar for all $R_{\min}{<}R{<}R_c$.  
We show the canonical path from the FCSS to the PCSS, where 
now the strategy to reach a patterned state 
already indicated in Fig\ref{vf2}(c) 
becomes more prominent.
Up to $t\approx 10$, $E$ is reduced around $x_2=-5$ and $x_4=5/3$. 
Conversely, $E$ is initially increased near the minima $x_1=-5/3$ and 
$x_3=5$ of $\hat v$, leading to a decrease of $v$ and an increase 
of $w$ near $x_{1,3}$. 
On the other hand, due to diffusion of $w$ 
the increase of $v$ near $x_{2,4}$ does {\em not} lead to a decrease 
of $w$ compared to the FCSS $w_0$. Instead $w$ increases significantly 
{\em everywhere}. After this transient 
the harvesting effort $E$ is increased near $x_{2,4}$, leading to 
an overall quick convergence of $u(t)$ to the PCSS $\uh$. 

Thus, the main point for the strategy to go from $\uh_{\text{FCSS}}$ 
to $\uh_{\text{PS}}$  
is to initially introduce a spatial variation (of the right 
wavelength) into $E$, which yields maxima of $v$ at the minima of this initial 
$E$, but then to rather quickly switch to the harvesting on the slopes 
of the generated maxima of $v$. The canonical path shows precisely how to do 
this. 
Also note (blue curve in (Fig.~\ref{f2b}a)) that the initial harvesting briefly yields a higher 
$J_{c,a}$  than $J_{c,a}(\uh_{\text{PS}})$ 
but in the transition $4<t<25$, say, $J_{c,a}(t)$ is significantly 
below $J_{c,a}(\uh_{\text{PS}})$. 

For the values we have 
\huga{\label{r10val}
J(\uh_{\text{FCSS}})=116.9<J(u_\text{FCSS$\to$PS})=132.49<
J(\uh_{\text{PS}})=145.26. 
}
Thus, again tentatively, 
see Remark \ref{prerem}, and in particular \reff{p0val} below, 
we classify $\uh_{\text{PS}}$ at $R=10$ as a POSS, with 
$\uh_{\text{FCSS}}$ in its domain of attraction. 
 For completeness we remark that at $R=20$ we have 
\huga{\label{r20val}
\text{$J(\uh_{\text{FCSS}})=455.31$, 
$J(\uh_{\text{PS}})=480.88$, and $J(u_{\text{FCSS$\to$PS}})=474.57$.}
}

\subsubsection{A patterned Skiba point}\label{pskibasec}
In ODE OC applications, if there are several locally stable OSS, 
then often an important issue is to identify their domains of attractions. 
These are separated by so called threshold or Skiba--points (if $N=1$) 
or Skiba--manifolds (if $N>1$) \cite{Skiba78, grassetal2008,KW10}. 
\begin{figure}[ht]
\vs{-32mm}
\bce \ig[width=0.99\textwidth, height=110mm]{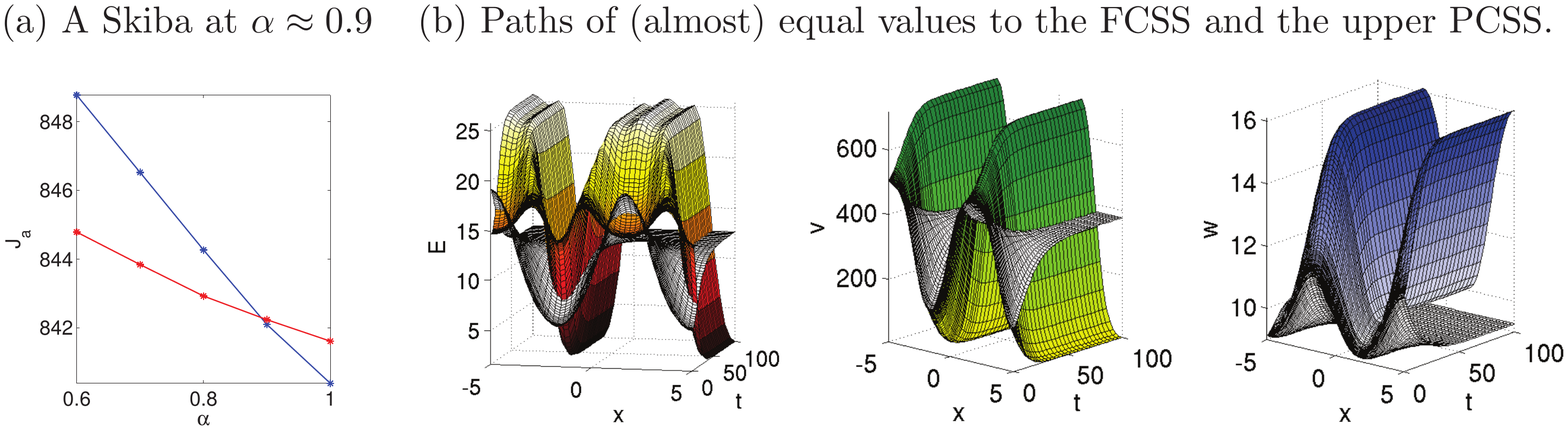}

\vs{-35mm}
   \caption{{\small In (a), the blue and red 
lines gives $J$ for the canonical paths $u_{\text{$\to$FCSS}}$ and 
$u_{\text{$\to$PS}}$ 
from $(v,w)_\al(0):=\al (v,w)_{\text{PS}}+(1-\al)(v,w)_{\text{FCSS}}$ 
to FCSS/pt13 and to p1/pt34, respectively. In (b), the white surfaces 
are for $u_{\text {$\to$FCSS}}$ and the colored ones for $u_{\text {$\to$PS}}$.
}\label{vf3}}
\ece
\end{figure}
Roughly speaking, these are (consist of) initial states from which 
there are more than one 
optimal paths {\em with the same value, but connecting to different CSS}. 
In PDE 
applications, even under spatial discretization with moderate $nN$, 
Skiba manifolds should be expected to become very complicated objects. 

In Fig.~\ref{vf3} we just given one example of a patterned Skiba 
point ``between'' $\uh_{\text{PS}}$ given by p1/pt34 
and $\uh_{\text{FSS}}$ given by 
FCSS/pt13, at $R=28$, where $\uh_{\text{PS}}$ and $\uh_{\text{FCSS}}$ are 
the two 
possible targets for canonical paths. The blue and red lines in (a) gives 
$J$ for the canonical paths as indicated.  
The lines intersect near $\al\approx 0.9$, 
giving the same value $J$. Hence, while the two paths are completely 
different, they both are equally optimal, and 
for illustration (b) shows the two paths for $\al=0.9$, where 
$|J_{\text{$\to$FCSS}}-J_{\text{$\to$PS}}|<0.08$.

\subsection{2D results}\label{2dsec}
The basic mechanisms of pattern formation in 
reaction--diffusion models can usually be studied in 1D, but 
for quantitative results for vegetation--water ecosystem models one should 
also consider the more pertinent 2D situation, and clearly the same 
applies to the OC system. Even though we do not 
expect qualitatively different results from the 1D case, here we 
give a short overview over 2D PCSS and the associated canonical paths.  

\begin{figure}[ht]\bce
\vs{-12mm}
\ig[width=0.99\textwidth, height=110mm]{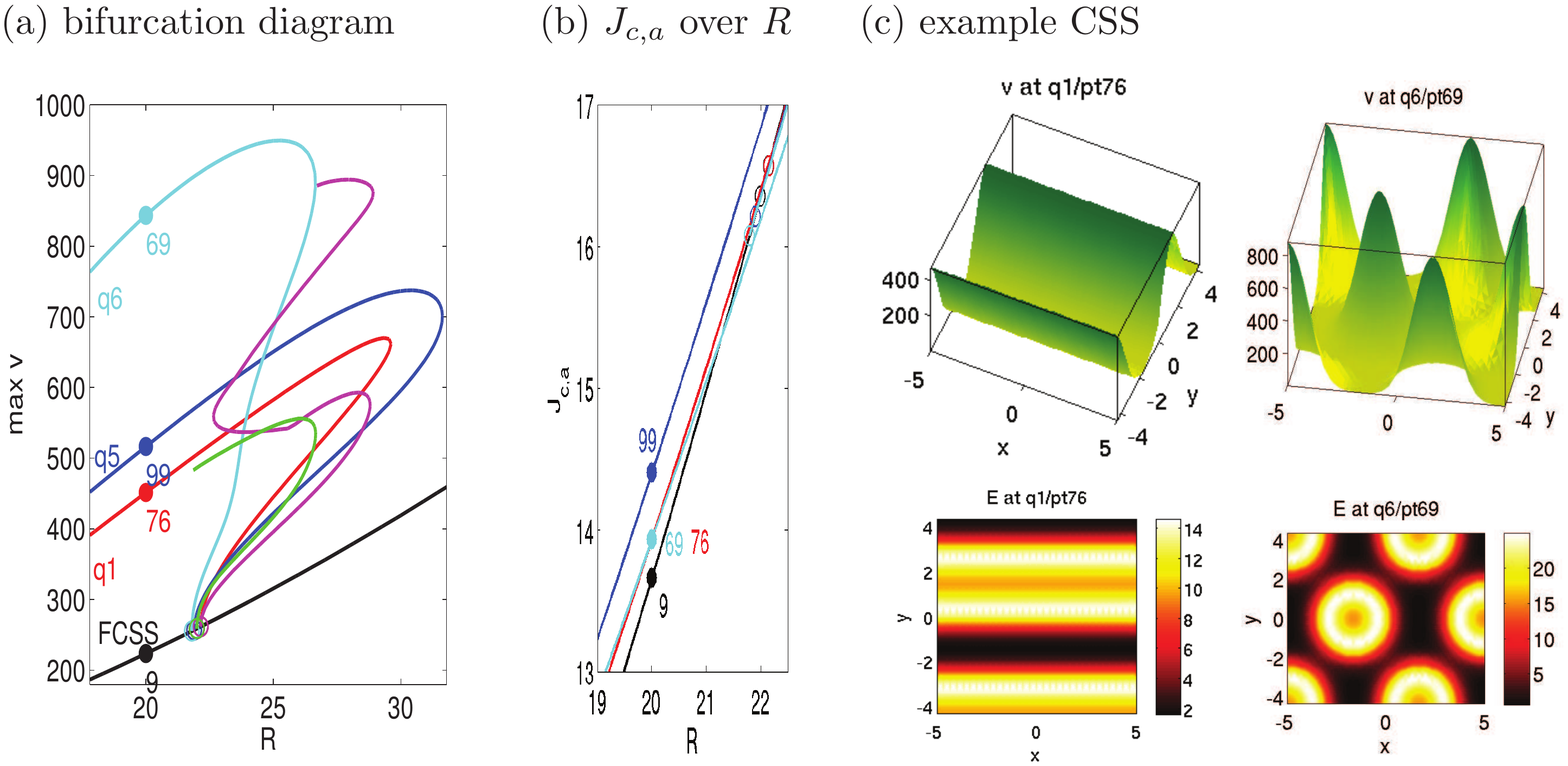}
\ece
\vs{-22mm}
\caption{{\small Partial bifurcation diagram and example plots of CSS in 2D.  
\label{f3}}}
\end{figure}

The first question concerns the second spatial length, now called 
$y$--direction. For classical 
Turing bifurcations, typically stripes and/or hexagonal spots are 
the most stable bifurcating 2D patterns, see \cite{uwsnak14} and the references 
therein for discussion, 
and thus by analogy here we choose the domain 
$\Om=(-L,L)\times (-\sqrt{3} L/2, \sqrt{3}L/2)$ with $L=5$ as in 1D. 
Figure \ref{f3}(a) shows five branches bifurcating from the FCSS. 
It turns out, that the 1D branch p1 actually comes out of the 5th 
bifurcation point in 2D, and is therefore called q5 now, while the 
first branch q1 corresponds to horizontal stripes, see (c). 
Thus, $L=5$, chosen for simplicity, does 
not capture the first instability in 1D. However, the first 
bifurcation points are very close together; moreover, upon continuation 
to low $R$ the q5 branch still yields the highest $J_{c,a}$, see (b).

The sixth branch q6 is a hexagon branch.  
Similar to q5, both q1 and q6 have the SPP after their first folds. 
The other branches are other types of stripes or spots, 
for instance ``squares'', but none of these have the SPP. All branches exhibit 
some secondary bifurcations, not shown here, 
and to not overload the bifurcation 
diagram we only plot the starting segments of q2 (green) and q3 (magenta). 

\begin{figure}
\vs{-30mm}
\bce \ig[width=0.7\textwidth, height=150mm]{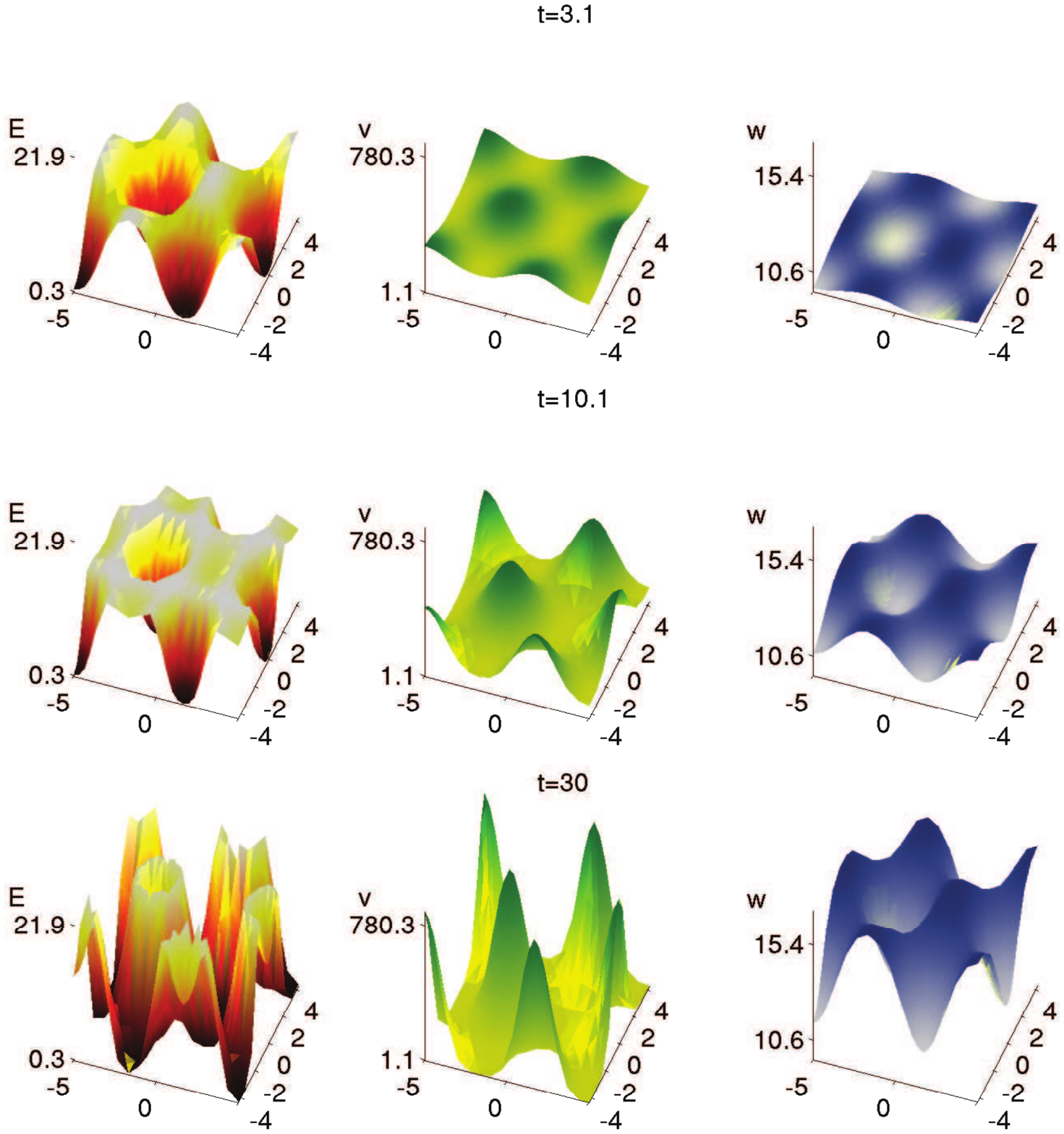}
\ece 
\vs{-30mm}
\caption{{\small Snapshots of $E$ and  $(v,w)$ on the canonical 
path from the FCSS to the 
hexagons. The states $v,w$ both go directly into their ultimate hex 
pattern, which then only grows, while $E$ shows a switching 
analogously to the 1D case.} \label{f4}}
\end{figure}

At,  for instance, $R=20$ we now have 3 possible targets for canonical 
paths: $\uh_{\text{hs}}$ (horizontal stripes) from q1, 
$\uh_{\text{hex}}$ (hexagons) from q6, and $\uh_{\text{vs}}$ 
(vertical stripes)  from q5, already 
discussed in 1D as p1. It turns out that we can reach each of these 
from the FCSS, with $J(\uh_{{\rm FCSS}})=455.31$, cf.~\reff{r20val}, 
and similarly $J(\uh_{\text{vs}})=480.88$ and $\uh_{\text{FCSS$\to$vs}}=474.56$  
 as in 1D (since these values are normalized by $|\Om|$). 
For the path to $\uh_{\text{hs}}$ with $J(\uh_{\text{hs}})=464.33$ we obtain 
$J=464.91$, and the path to $\uh_{\text{hex}}$ with $J(\uh_{\text{hex}})
=464.33=J(\uh_{\text{hs}})$ 
(up to 2 digits) yields $J=467.38$. Thus we again have $V(\uh_{{\rm FCSS}})=
474.6$. The strategies for these paths are the natural extensions 
of the 1D case: given a target PCSS $\uh$, initially $E$ has 
minima at the maxima of $\hat v$, but after a rather short transient 
during which $v(t,\cdot)$ develops maxima at the right places, 
$E$ changes to harvesting in the neighborhood of these maxima. 
Movies of these paths can be found at \cite{p2phome}, 
and in Fig.~\ref{f4} 
we present some snapshots. 
We could not compute canonical paths from any of the PCSS to any 
other PCSS, with the continuation typically failing due to a sequence of 
folds. 
Thus we strongly expect all three PCSS to be locally stable POSS.

\subsection{Remarks on further parameter dependence}
\label{rhosec}
So far we varied the rainfall $R$ as our external bifurcation parameter. 
Similarly, we could vary some other of the physical parameters 
$g,\eta,\ldots,d_{1,2}$ 
(first six rows of Table \ref{tab1}), and in most cases may expect 
bifurcations to patterned CSS. 

Maybe even more interesting from an application point of view 
is the dependence on the economic parameters 
$\rho, c,p$ and $\al$ (discount rate, cost for harvesting, price of harvest, 
and elasticity), as these may vary strongly with economic 
circumstances. 
Moreover, varying a second parameter often 
also gives bifurcations 
to branches which were missed upon continuation of just the primary 
parameter, and these may play an important role in the overall 
structure of the solution set; this does happen here, see \S\ref{unisec} below.

\subsubsection{Experiments with the discount rate $\rho$}\label{dratesec}
\begin{figure}[ht]
\bce
\vs{-12mm}
\ig[width=0.8\textwidth, height=95mm]{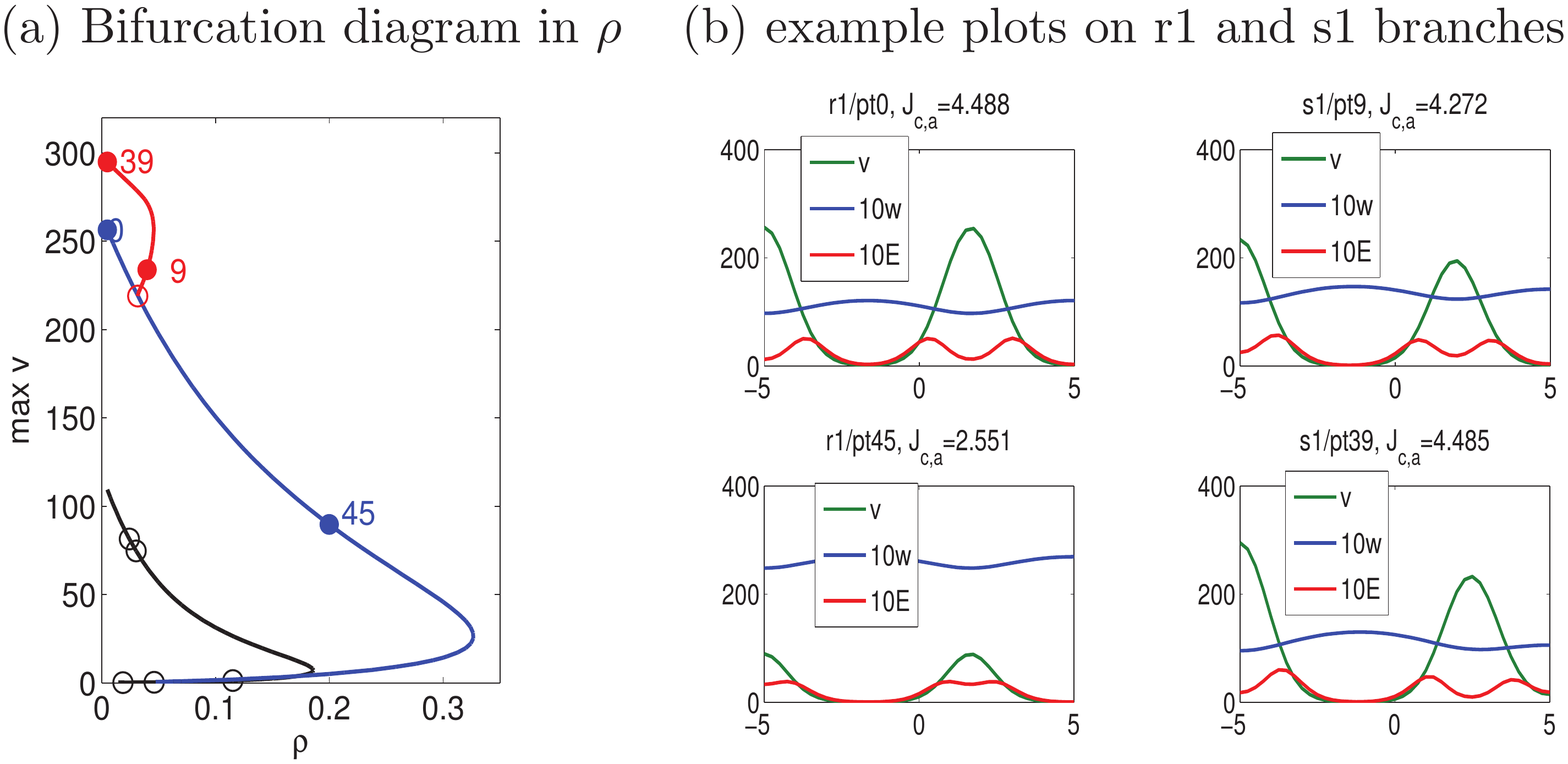}
\ece
\vs{-25mm}
   \caption{{\small Bifurcations when varying the 
discount rate $\rho$ at $R=10$. 
The blue branch (which at $\rho=0.03$ corresponds to Fig.~\ref{f2b}) 
loses the SPP at $\rho\approx 0.032$, where the red s1-branch bifurcates, 
which has the SPP up to its fold. }\label{rf1}}
\end{figure}

In Fig.\ref{rf1} we illustrate the dependence  of 
the PCSS on the p1 branch from Fig.~\ref{f1} on $\rho$, at fixed $R=10$, 
cf.~also Fig.~\ref{f2b}. Panel (a) shows the bifurcation diagram; 
to obtain the blue and black branches we continued 
the points p1/pt65 and FCSS/pt44 from $\rho=0.03$ down to $\rho=0.005$, 
reset counters, renamed p1 to r1, and continued to larger $\rho$ again. 
Both, the FCSS and the r1 branches then show some folds at $\rho\approx 0.185$
and $\rho\approx 0.325$, respectively. 
More importantly, the r1 branch has the SPP 
for small $\rho$, but loses it 
at $\rho\approx 0.032$ to another PCSS branch s1. Solutions on s1 
have maxima of different heights in $v$, see (b), and have the SPP only 
up to the 
fold at $\rho=\rho_f\approx 0.046$. Moreover, there are further bifurcations 
from the FCSS branch to PCSS branches, but none of these has the SPP. 

With $\uh_{\text{PS}}$ given by s1/pt9 ($\rho=0.04$), the control 
to go from the FCSS to $\uh$ shows a similar switching strategy 
as, e.g., in Fig.~\ref{f2b}. The values are 
\huga{\label{skewval}
J(\uh_{\text{FCSS}})=82.5,\quad J(u_{\text{FCSS$\to$PS}})=92.8,\quad 
J(\uh_{\text{PS}})=106.8,
}
which are more or less comparable 
to Fig.~\ref{f2b} (with $\rho=0.03$). 
On the other hand, for the canonical path from the FCSS to 
r1/pt0 at $\rho=0.005$ we obtain 
\huga{\label{lrval}
J(\uh_{\text{FCSS}})=774.9,\quad J(u_{\text{FCSS$\to$PS}})=892.2,
\quad J(\uh_{\text{PS}})=897.7.
} 
Additional to the larger 
total values due to the smaller discount rate, compared 
to Fig.~\ref{f2b} the canonical path to $\uh_{\text{PS}}$ now has 
almost the same value as $\uh_{\text{PS}}$ itself. 
This illustrates that at low $\rho$ the 
transients have less influence, as expected.  

For $\rho>\rho_f$ none of the CSS on 
the branches that are shown in Fig.~\ref{rf1}, or that can be obtained 
from the shown bifurcation points, have the SPP. This does not mean 
that PCSS with the SPP do not exist in this parameter regime, but rather 
that they must be obtained by continuation and bifurcation in some other way, 
cf.~Remark \ref{unirem} and \S\ref{unisec}. 

\subsubsection{Dependence on the price $p$, and the unimodal branch}\label{unisec}
\begin{figure}[ht]
\vs{-25mm}
\bce
\ig[width=0.99\textwidth, height=115mm]{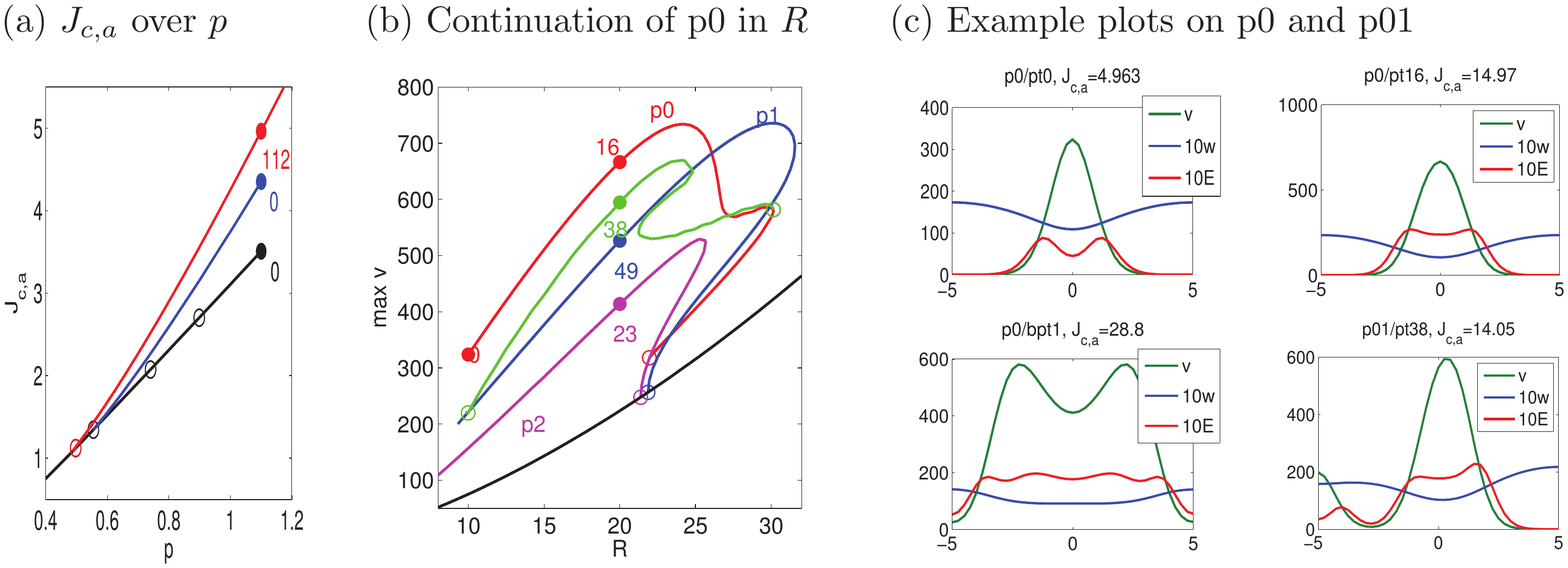}
\ece
\vs{-35mm}
   \caption{{\small (a) branches FCSS (black), p1 (blue), and p0 (red) over $p$, 
$R=10$ fixed. (b) Continuation in $R$ of new branch p0 found in (a), 
and bifurcating branch p01 (green), together with known branches FSS (black) and p2 from Fig.~\ref{f1}. Example plots in (c). 
}\label{rf2}}
\end{figure}

In Fig.~\ref{rf2}(a) we illustrate the dependence 
of the FCSS and p1 branches on the price $p$, with fixed $R=10$, 
starting at $p=1.1$ with p1/pt65 and FCSS/pt44 from Fig.~\ref{f1}, 
respectively.  
Naturally, the values 
decrease  as $p$ decreases, and not surprisingly p1 bifurcates 
from the FSS branch at some $p_c{\approx}0.55$. Next, as an additional 
benefit we find the ``unimodal'' branch p0, which bifurcates from the 
FCSS branch near $p{=}0.5$, and which yields a higher 
$J_{c,a}$ than p1. 

Thus, in a next step we continue p0 to $p=1.1$ and then 
switch back to continuation in $R$, with $p=1.1$ fixed, i.e., 
p0/pt0 in (b)-(c) is pt112 from (a). 
It turns out that the p0 branch has the SPP up the fold at 
$R_{\text{f}}\approx 30.15$, and slightly below the fold there is a 
bifurcation point to the green branch. This contains 
some ``skewed'' solutions, and connects p0 and p1. 
Ultimately, p0 connects back to p2 from 
Fig.~\ref{f1} at low amplitude near $R=21.1$. Thus, we could also have found 
p0 by following secondary bifurcations in Fig.~\ref{f1}. 
The values pertinent for 
the canonical path from FCSS/pt44 to 
p0/pt0 are 
\huga{\label{p0val}
J(\uh_{\text{FCSS}})=116.9\ (\text{as in \reff{r10val}}),\quad 
J(u_{\text{FCSS$\to$P0}})=145.11,\quad 
J(\uh_{\text{P0}})=165.42, 
}
which shows that the path to p0 dominates the path to p1 from 
Fig.~\ref{f2b}. Thus, the point FCSS/pt44 is in the domain of attraction 
of p0, and not of p1. 

On the other hand, we could not find canonical paths from p1/pt65 
to p0/pt0 (or vice versa). Therefore p1/pt65 can still be classified 
as an at least locally stable POSS, and similarly, p0/pt0 is {\em only  
locally stable}, since it does not attract p1/pt65. Next one could 
compute a number of Skiba--points (cf.~\S\ref{pskibasec}) 
``between'' p0 and p1 to roughly 
characterize the respective domains of attraction, but here we refrain 
from this. 

Finally, despite trying some further combinations 
of continuations/bifurcations and also suitable direct 
initial guesses followed by a Newton loop, 
we could not find a ``half''--modal branch in 
the parameter regions studied so far, 
i.e., a branch on which solutions are monotonous in $v$. Thus, it appears 
that such ``long wave'' PCSS, i.e., with period 20, 
do not exist in these parameter regions. 

\section{Summary and discussion}\label{dsec}
Our numerical approach for spatially distributed optimal control problems 
may yield rich results, if applied carefully, in the following sense. 
First, the canonical system may have many steady states, and 
it is in general not clear how to find all {\em relevant} CSS, 
and which of the CSS have the SPP and hence 
are suitable targets for canonical paths, 
and second 
it needs to be checked which CSS ultimately belong to optimal paths. 

 On the other hand, the value $J$ of the CSS 
can easily be calculated in parallel with the 
bifurcation diagram of the CSS, allowing to identify which 
CSS maximize $J$ amongst all CSS. Compared to the computation 
of canonical paths (or direct methods for the optimization 
problem \reff{oc1}), this first step is relatively cheap numerically, but  
(together with the SPP) 
typically still gives a strong indication for optimal CSS. 

The computation of canonical paths is a connecting orbits problem, and 
in particular in two spatial dimensions this may 
become numerically expensive. In practice we found our two--step 
approach to be reasonably fast for up to 5000 degrees 
of freedom of $u$ at fixed time, e.g., for our vegetation 
model 1250 spatial discretization 
points and 4 components, and 
up to 100 temporal discretization points; 
up to such values a continuation step in the calculation 
of a canonical path takes up to a minute on a desktop computer 
(Intel i7, 2.3 GHz), such that a typical canonical path is 
computed in about 5 continuation steps in at most 5 minutes, 
and often much more quickly. In 1D, with $n=50$, say, 
typical canonical paths are computed in 
a few seconds. 

Here we applied our method exemplarily to the optimal 
control model from \cite{BX10}, and for this we again 
summarize our main results as follows, cf.~also {\bf (a)--(d)} 
in the Introduction: 
\bci
\item[(i)] Compared to the case of private optimization, we have CSS (both flat 
and patterned) for significantly 
lower rainfall values $R$, and the whole Turing (like) bifurcation scenario 
is shifted to lower $R$. This is important for welfare as it means a much 
increased robustness of vegetation (and hence harvest) with respect to low 
rainfall.  
Moreover, at a given $R$ the social control gives a significantly 
(almost ten times) higher $J$ 
for steady states than private optimization; see, e.g., Table \ref{tab3}, bottom center, and Fig.~\ref{fnc}(c), for numerical values. 
\item[(ii)] At low $R$, some PCSS yield a higher $J$ than the FCSS, and 
some of these PCSS are locally stable POSS. 
\item[(iii)]  
The optimal controls to reach such a POSS $\uh$ from a FCSS 
follow some general rules: first decrease the harvesting effort $E$ at the 
location of the maxima of the desired POSS, then start harvesting 
near (but not at) the maxima of the POSS, as determined by $E$ from the POSS. 
The increase in welfare by controlling the system from a FCSS 
to a POSS can be up to 40\%, see \reff{p0val}. See also, e.g., \reff{cpv2}, \reff{r10val}, \reff{r20val}, and the values given in \S\ref{2dsec} for 
further numerical values. 
{\item[(iv)] The co--state (or shadow price) $\lam$ computed for an optimal path 
can be interpreted 
as an optimal tax for private optimization, see Remark \ref{taxrem}. }
\eci
{The points (ii) and (iii) emphasize that resource management rules 
in systems with multiple CSS may need to take several of these into account, 
as also illustrated by the Skiba point computations in \S\ref{pskibasec}.} 
We strongly expect similar results about POSS in other spatially 
distributed optimal control problems for (Turing like) systems 
of PDE. Thus we hope  that our numerical approach 
is a valuable tool to 
study the basic behavior of spatially distributed optimal control models  
with the states fulfilling a  reaction diffusion system.

\renewcommand{\arraystretch}{1.1}\renewcommand{\baselinestretch}{1.05}
\small
\input{veg-r1.bbl}

\end{document}

%% file: mindef.tex
\def\bmip{\begin{minipage}{\textwidth}}\def\emip{\end{minipage}}
\def\huga#1{\begin{gather} #1 \end{gather}}

\def\hual#1{\begin{align} #1 \end{align}}

\newcommand{\R}{{\mathbb R}}

\newcommand{\argmax}{\operatornamewithlimits{argmax}}

\def\ig{\includegraphics}

\def\CD{{\cal D}}  
\def\CH{{\cal H}}

\def\del{\delta}
\def\ga{\gamma}\def\uh{\hat{u}}
\def\noi{\noindent}\def\ds{\displaystyle}\def\al{\alpha}
\def\pa{{\partial}}\def\lam{\lambda}
\def\Om{\Omega}
\def\dd{\, {\rm d}}\def\er{{\rm e}}
\newcommand{\bi}{\begin{itemize}}\newcommand{\ei}{\end{itemize}}
\newcommand{\ben}{\begin{enumerate}}\newcommand{\een}{\end{enumerate}}
\newcommand{\bci}{\begin{compactitem}}\newcommand{\eci}{\end{compactitem}}
\newcommand{\bcen}{\begin{compactenum}}\newcommand{\ecen}{\end{compactenum}}
\newcommand{\bce}{\begin{center}}\newcommand{\ece}{\end{center}}
\newcommand{\reff}[1]{(\ref{#1})}
\newcommand{\ov}[1]{{\overline {#1}}}
\newcommand{\spr}[1]{\left\langle #1 \right\rangle}
\newcommand{\hs}[1]{{\hspace{#1}}}\newcommand{\vs}[1]{{\vspace{#1}}}

\def\ra{\rightarrow}

\newcommand{\barr}{\begin{array}}\newcommand{\earr}{\end{array}}
\newcommand{\bpm}{\begin{pmatrix}}\newcommand{\epm}{\end{pmatrix}}
\newcommand{\bsm}{\left(\begin{smallmatrix}}
\newcommand{\esm}{\end{smallmatrix}\right)}
\newcommand{\ba}{\begin{array}}\newcommand{\ea}{\end{array}}

%% file: veg-r1.bbl
\newcommand{\etalchar}[1]{$^{#1}$}

%% file: veg-r1.bbl
\begin{thebibliography}{RDdRvdK04}

\bibitem[AAC11]{AAC11}
S.~Ani{\c{t}}a, V.~Arn{\u{a}}utu, and V.~Capasso.
\newblock {\em An introduction to optimal control problems in life sciences and
  economics}.
\newblock Birkh\"auser/Springer, New York, 2011.

\bibitem[ACKLT13]{ACKT13}
S.~Ani{\c{t}}a, V.~Capasso, H.~Kunze, and D.~La~Torre.
\newblock Optimal control and long-run dynamics for a spatial economic growth
  model with physical capital accumulation and pollution diffusion.
\newblock {\em Appl. Math. Lett.}, 26(8):908--912, 2013.

\bibitem[ADS14]{Apre14}
N.~Apreutesei, G.~Dimitriu, and R.~Strugariu.
\newblock An optimal control problem for a two-prey and one-predator model with
  diffusion.
\newblock {\em Comput. Math. Appl.}, 67(12):2127--2143, 2014.

\bibitem[BEGX13]{Brocketal2013}
W.A. Brock, G.~Engstr\"om, D.~Grass, and A.~Xepapadeas.
\newblock Energy balance climate models and general equilibrium optimal
  mitigation policies.
\newblock {\em Journal of Economic Dynamics and Control}, 37(12):2371--2396,
  2013.

\bibitem[BPS01]{BPS01}
W.J. Beyn, Th. Pampel, and W.~Semmler.
\newblock Dynamic optimization and {S}kiba sets in economic examples.
\newblock {\em Optimal Control Applications and Methods}, 22(5--6):251--280,
  2001.

\bibitem[BX08]{BX08}
W.A. Brock and A.~Xepapadeas.
\newblock Diffusion-induced instability and pattern formation in infinite
  horizon recursive optimal control.
\newblock {\em Journal of Economic Dynamics and Control}, 32(9):2745--2787,
  2008.

\bibitem[BX10]{BX10}
W.~Brock and A.~Xepapadeas.
\newblock Pattern formation, spatial externalities and regulation in coupled
  economic--ecological systems.
\newblock {\em Journal of Environmental Economics and Management},
  59(2):149--164, 2010.

\bibitem[Cla90]{Cl90}
C.~W. Clark.
\newblock {\em Mathematical bioeconomics}.
\newblock John Wiley \& Sons, Inc., New York, second edition, 1990.
\newblock The optimal management of renewable resources.

\bibitem[CPB12]{CP12}
C.~Camacho and A.~P\'erez-Barahona.
\newblock {Land use dynamics and the environment. Documents de travail du
  Centre d’Economie de la Sorbonne}, 2012.

\bibitem[DBC{\etalchar{+}}08]{DBC08}
V.~Deblauwe, N.~Barbier, P.~Couteron, P.~Lejeune, and J.~Bogaert.
\newblock The global biogeography of semi-arid periodic vegetation patterns.
\newblock {\em Global Ecol. Biogeogr.}, 17:715--723, 2008.

\bibitem[DL09]{DL09}
W.~Ding and S.~Lenhart.
\newblock Optimal harvesting of a spatially explicit fishery model.
\newblock {\em Natural Resource Modeling}, 22(2):173--211, 2009.

\bibitem[DRUW14]{p2p2}
T.~Dohnal, J.~Rademacher, H.~Uecker, and D.~Wetzel.
\newblock {pde2path 2.0}.
\newblock In H.~Ecker, A.~Steindl, and S.~Jakubek, editors, {\em {ENOC 2014 -
  Proceedings of 8th European Nonlinear Dynamics Conference, ISBN:
  978-3-200-03433-4}}, 2014.

\bibitem[GCF{\etalchar{+}}08]{grassetal2008}
D.~Grass, J.P. Caulkins, G.~Feichtinger, G.~Tragler, and D.A. Behrens.
\newblock {\em Optimal Control of Nonlinear Processes: With Applications in
  Drugs, Corruption, and Terror}.
\newblock Springer Verlag, 2008.

\bibitem[Gra15]{grass2014}
D.~Grass.
\newblock From {0D} to {1D} spatial models using {OCMat}.
\newblock Technical report, ORCOS, 2015.

\bibitem[GRS14]{GRS14}
K.~Gowda, H.~Riecke, and M.~Silber.
\newblock Transitions between patterned states in vegetation models for
  semiarid ecosystems.
\newblock {\em PRE}, 89:022701, 2014.

\bibitem[GU15]{GU15}
D.~Grass and H.~Uecker.
\newblock Optimal management and spatial patterns in a distributed shallow lake
  model.
\newblock {P}reprint, 2015.

\bibitem[HRvdB{\etalchar{+}}01]{riet01}
R.~HillerisLambers, M.G. Rietkerk, F.~van~den Bosch, H.H.T. Prins, and
  H.~de~Kroon.
\newblock Vegetation pattern formation in semi-arid grazing systems.
\newblock {\em Ecology}, 82:50--61, 2001.

\bibitem[KW10]{KW10}
T.~Kiseleva and F.O.O. Wagener.
\newblock Bifurcations of optimal vector fields in the shallow lake system.
\newblock {\em Journal of Economic Dynamics and Control}, 34(5):825--843, 2010.

\bibitem[KXL15]{KXL15}
M.~R. Kelley, Y.~Xing, and S.~M. Lenhart.
\newblock Optimal fish harvesting for a population modeled by a nonlinear
  parabolic partial differential equation.
\newblock {\em Natural Resource Modeling, DOI 10.1111/nrm.12073}, 2015.

\bibitem[LM01]{LM01}
S.~M. Lenhart and J.~A. Montero.
\newblock Optimal control of harvesting in a parabolic system modeling two
  subpopulations.
\newblock {\em Math. Models Methods Appl. Sci.}, 11(7):1129–1141, 2001.

\bibitem[LW07]{LW07}
S.~Lenhart and J.~Workman.
\newblock {\em Optimal Control Applied to Biological Models}.
\newblock Chapman Hall, 2007.

\bibitem[MT04]{MT04}
F.~Mazzia and D.~Trigiante.
\newblock A hybrid mesh selection strategy based on conditioning for boundary
  value {ODE} problems.
\newblock {\em Numerical Algorithms}, 36(2):169--187, 2004.

\bibitem[Mur89]{Murray89book}
J.~D. Murray.
\newblock {\em Mathematical {B}iology}.
\newblock Springer, Berlin, 1989.

\bibitem[Neu03]{Neu03}
M.~G. Neubert.
\newblock Marine reserves and optimal harvesting.
\newblock {\em Ecology Letters}, 6(9):843--849, 2003.

\bibitem[NPS11]{NPS11}
I.~Neitzel, U.~Pr{\"u}fert, and Th. Slawig.
\newblock A smooth regularization of the projection formula for constrained
  parabolic optimal control problems.
\newblock {\em Numer. Funct. Anal. Optim.}, 32(12):1283--1315, 2011.

\bibitem[QB12]{QB12}
M.F. Quaas and S.~Baumgärtner.
\newblock Optimal grazing management rules in semi--arid rangelands with
  uncertain rainfall.
\newblock {\em Natural Resource Modeling}, 25(2):364--387, 2012.

\bibitem[RDdRvdK04]{riet04}
M.G. Rietkerk, S.C. Dekker, P.C. de~Ruiter, and J.~van~de Koppel.
\newblock Self-organized patchiness and catastrophic shifts in ecosystems.
\newblock {\em Science}, 305:1926--1929, 2004.

\bibitem[RZ99a]{RZ99}
J.~P. Raymond and H.~Zidani.
\newblock Hamiltonian {P}ontryagin's principles for control problems governed
  by semilinear parabolic equations.
\newblock {\em Appl. Math. Optim.}, 39(2):143--177, 1999.

\bibitem[RZ99b]{RZ99b}
J.~P. Raymond and H.~Zidani.
\newblock Pontryagin's principle for time-optimal problems.
\newblock {\em J. Optim. Theory Appl.}, 101(2):375--402, 1999.

\bibitem[SBB{\etalchar{+}}09]{SBB09}
M.~Scheffer, J.~Bascompte, W.~A. Brock, V.~Brovkin, St.~R. Carpenter, V.~Dakos,
  H.~Held, E.~H. van Nes, M.~Rietkerk, and G.~Sugihara.
\newblock Early-warning signals for critical transitions.
\newblock {\em Nature}, 461:53--59, 2009.

\bibitem[Ski78]{Skiba78}
A.~K. Skiba.
\newblock Optimal growth with a convex-concave production function.
\newblock {\em Econometrica}, 46(3):527--539, 1978.

\bibitem[SZvHM01]{meron01}
M.~Shachak, Y.~Zarmi, J.~von Hardenberg, and E.~Meron.
\newblock Diversity of vegetation patterns and desertification.
\newblock {\em PRL}, 87, 2001.

\bibitem[Tr{\"o}10]{Tr10}
Fredi Tr{\"o}ltzsch.
\newblock {\em Optimal control of partial differential equations}, volume 112
  of {\em Graduate Studies in Mathematics}.
\newblock American Mathematical Society, Providence, RI, 2010.

\bibitem[Uec14]{p2phome}
H.~Uecker.
\newblock {pde2path, {\tt www.staff.uni-oldenburg.de/hannes.uecker/pde2path}},
  2014.

\bibitem[Uec15]{U15p2}
H.~Uecker.
\newblock The pde2path add-on library p2poc for solving o infinite
  time–horizon spatially distributed optimal control problems - {Quickstart
  Guide}.
\newblock {P}reprint, 2015.

\bibitem[UW14]{uwsnak14}
H.~Uecker and D.~Wetzel.
\newblock {Numerical results for snaking of patterns over patterns in some 2D
  Selkov-Schnakenberg Reaction-Diffusion systems}.
\newblock {\em SIADS}, 13-1:94--128, 2014.

\bibitem[UWR14]{p2pure}
H.~Uecker, D.~Wetzel, and J.~Rademacher.
\newblock {pde2path -- a Matlab package for continuation and bifurcation in 2D
  elliptic systems}.
\newblock {\em {NMTMA}}, 7:58--106, 2014.

\bibitem[Xep10]{Xe10}
A.~Xepapadeas.
\newblock The spatial dimension in environmental and resource economics.
\newblock {\em Environment and Development Economics,
  doi:10.1017/S1355770X10000355}, 2010.

\bibitem[ZKY{\etalchar{+}}13]{meron13}
Y.~Zelnik, S.~Kinast, H.~Yizhaq, G.~Bel, and E.~Meron.
\newblock Regime shifts in models of dryland vegetation.
\newblock {\em Philos. Trans. R. Soc. Lond. Ser. A Math. Phys. Eng. Sci.},
  371(2004), 2013.

\end{thebibliography}
